\documentclass[12pt]{article}
\usepackage{setspace}
\usepackage[utf8]{inputenc}
\usepackage{amsmath,amsfonts}
\usepackage{color,enumerate}
\usepackage{graphics,afterpage,subfigure,graphicx,subfigure,lipsum,caption,enumerate}
\usepackage[utf8]{inputenc}
\usepackage{mathtools}   
\usepackage{pifont}
\usepackage[round]{natbib} 
\usepackage[OT1]{fontenc}
\usepackage{amscd}
\usepackage{latexsym,amssymb,amsthm}
\usepackage{pdflscape}
  \usepackage{hyperref}
  \hypersetup{colorlinks=true,citecolor=blue, backref,linktocpage=true,linkcolor=blue}

\newcommand{\yn}{{\bf y}_n}

\newcommand{\pg}{p_{\bgm}}
\newcommand{\bgm}{\boldsymbol{\gamma}}
\newcommand{\bbta}{\boldsymbol{\beta}}
\newcommand{\Rg}{R_{\bgm}}
\newcommand{\ad}{\alpha_\delta} 
\newcommand{\bta}{\boldsymbol{\theta}}

\newtheorem{lemma}{Lemma}
\newtheorem{remark}{Remark}
\newtheorem{theorem}{Theorem}

\title{\bf Targeted Random Projection for Prediction from High-Dimensional Features}
\author{Minerva Mukhopadhyay \and David B. Dunson}
\date{}

\begin{document}

\maketitle
\begin{abstract}
We consider the problem of computationally-efficient prediction from high dimensional and highly correlated predictors in challenging settings where accurate variable selection is effectively impossible. Direct application of penalization or Bayesian methods implemented with Markov chain Monte Carlo can be computationally daunting and unstable. Hence, some type of dimensionality reduction prior to statistical analysis is in order. Common solutions include application of screening algorithms to reduce the regressors, or dimension reduction using projections of the design matrix. The former approach can be highly sensitive to threshold choice in finite samples, while the later can have poor performance in very high-dimensional settings. We propose a TArgeted Random Projection (TARP) approach that combines positive aspects of both strategies to boost performance. In particular, we propose to use information from independent screening to order the inclusion probabilities of the features in the projection matrix used for dimension reduction, leading to data-informed sparsity. We provide theoretical support for a Bayesian predictive algorithm based on TARP, including both statistical and computational complexity guarantees. Examples for simulated and real data applications illustrate gains relative to a variety of competitors. 
\end{abstract}

\par \noindent
{\bf Some key words:} Bayesian; Dimension reduction; High-dimensional; Large $p$, small $n$; Random projection; Screening.
\par\medskip\noindent
{\bf Short title}: Targeted Random Projection


\section{Introduction}\label{sec:1}

In many applications, the focus is on prediction of a response variable $y$ given a massive-dimensional vector of predictors ${\bf x}=(x_1,x_2, \ldots x_p)^{\prime}$. Often enormous numbers of possibly collinear predictors ${\bf x}$ are collected, and the sample size $n$ is modest relative to $p$, so that $p\gg n$. In such situations, it is common to assume that ${\bf x}$ can be replaced by a much lower-dimensional feature vector comprised of sparse linear combinations of the original predictors.  However, accurate learning of the precise lower-dimensional structure is often not possible, as the data simply do not contain sufficient information even putting aside the intractable computational problem.

There is a large literature on variable selection in $p \gg n$ settings.  Much of the focus has been on penalized optimization-based approaches, with some popular methods including LASSO \citep{lasso}, SCAD \citep{SCAD}), elastic net \citep{ENet}, the Dantzig selector \citep{DS}, and MCP \citep{MCP}. There is also a large literature on Bayesian approaches that attempt to characterize uncertainty in variable selection. Most approaches use some variation on the spike and slab prior (e.g. \cite{SpikeSlab}) or continuous shrinkage priors concentrated at zero with heavy tails (e.g. \cite{horseshoe}).  There is an increasingly rich theoretical literature providing asymptotic support for such methods in the $p \to \infty$ as $n \to \infty$ setting.  However, positive results rely on strong assumptions in terms of a high signal-to-noise ratio and low linear dependence in the columns of the design matrix.  

We are interested in settings where practical performance of the above methods is poor, due to a combination of statistical and computational intractability.  In such settings, it is common to use variable screening as a pre-processing step. In particular, independent screening tests for association between $y$ and each $x_j$ separately, and selects predictors with the largest or most significant associations for second stage analysis.  In general, screening can be guaranteed {\em asymptotically} to select a superset of the `true' predictors \citep{fan2009}.  When the number of predictors is sufficiently reduced, one can apply a simple maximum likelihood approach, penalized optimization, or Bayesian Markov chain Monte Carlo (MCMC) algorithms in the second stage.  However, when the predictors are highly correlated and/or the true data generating process does not exhibit strong sparsity with a high signal-to-noise ratio, it may be necessary to use a very conservative threshold for the measure of marginal association, limiting the dimensionality reduction occurring in the first stage.

As an alternative to variable screening, there is a rich literature on using random projections (RPs) to reduce data dimensionality prior to statistical analysis.  For example, {\em compressive sensing} uses RPs to reduce storage and communication costs in signal processing.  By exploiting sparsity, the original signal can be reconstructed from the compressive measurements with high accuracy (see, e.g., \cite{CT2005}, \cite{donoho2006}, \cite{davenport2010}). Usual compressive sensing acts in a row-wise manner, reducing the dimensionality of the design matrix from $n \times p$ to $m \times p$, with $m\ll n$.  This does not solve the big $p$ problem.  There is a relatively smaller literature on column-wise compression, which instead reduces the design matrix from $n \times p$ to $n \times m$, with $m \ll p$, while providing bounds on predictive errors (see, e.g., \cite{maillard2009}, \cite{fard2012}, \cite{kaban2014}, \cite{thanei2017}, \cite{GD_2015}, \cite{koop2016}). \cite{GD_2015} concentrate on approximating predictive distributions in Bayesian regression. The above literature on RPs focuses primarily on random matrices with i.i.d elements.

When predictors are very high-dimensional, existing RP methods can fail as they tend to include many unimportant predictors in each linear combination, diluting the signal. Potentially, one can attempt to improve performance by estimating the projection matrix, but this results in a daunting computational and statistical problem. Alternatively, we propose a TArgeted Random Projection (TARP) approach, which includes predictors in the RP matrix with probability proportional to their marginal utilities. These utilities are estimated quickly in a first stage using an independent screening-type approach.  To reduce sensitivity of the results to the different realizations of the RP matrices, we aggregate over multiple realizations.  TARP can be viewed as a type of rapid preconditioning, enabling improved predictive performance in high-dimensional settings. Compared with applying RPs after screening out predictors, TARP has the advantage of removing sensitivity to threshold choice by using a soft probabilistic approach.  

In Section \ref{sec:2}, we propose the methodology including the computational algorithm, choice of different tuning parameters and an analysis of computational complexity. 
We focus on generalized linear models (GLMs) for ease in presentation and development of a strong theoretical justification, although TARP can be applied directly in general settings to provide lower-dimensional features that can then be used in any predictive algorithm (random forests, Gaussian processes, neural nets, etc).  Section \ref{sec:3} provides theory on convergence rates for the predictive distribution of $y$.  Section \ref{sec:4} contains a simulation study comparing performance with a variety of competitors.  In Section \ref{sec:5}, we apply TARP to a variety of real data applications including a genomics dataset with millions of predictors. Section \ref{sec:7} contains a brief discussion, and proofs are included in an Appendix.

\section{The Proposed Method} \label{sec:2}
Let $\mathcal{D}^n=\left\{(\yn; X_n): \yn \in \mathbb{R}^n, X_n \in \mathbb{R}^{n\times p_n} \right\}$ denote the dataset consisting of $n$ observations on $p_n$ predictors $x_1, x_2, \ldots, x_{p_n}$ and a response $y$, and $(y_i; {\bf x}_i )$ denote the $i^{th}$ data point, $i=1,2, \ldots, n$. Suppose that the data can be characterized by a generalized linear model (GLM). The density of $y$ is related to the predictors as 
\begin{equation}
f(y_i| \bbta, \sigma^2)=\exp\left[\frac{1}{d(\sigma^2)} \left\{ y_i a({\bf x}_i^{\prime} \bbta) + b({\bf x}_i^{\prime} \bbta) +c(y_i) \right\} \right], \label{eq_rp1}
\end{equation}
where $a(\cdot)$ and $b(\cdot)$ are continuously differentiable functions, $a(\cdot)$ has non-zero derivative and $d(\cdot)$ is a non-zero function. The vector of regression coefficients $\bbta\in \mathbb{R}^{p_n}$, and $\sigma^2\in \mathbb{R}^{+}$ is the scale parameter. We approximate the density of $y$ in a compressed regression framework as follows: 
\begin{equation*}
 f(y_i | \bta, R_n, \sigma^2 )=\exp \left[ \frac{1}{d(\sigma^2)} \left\{ y_i a \left((R_n {\bf x}_i)^{\prime} \bta \right) + b\left((R_n {\bf x}_i)^{\prime} \bta\right) +c(y_i) \right\}\right].
\end{equation*}
Here $R_n \in \mathbb{R}^{m_n \times p_n}$ is a random projection matrix, $\bta\in \mathbb{R}^{m_n}$ is the vector of compressed regression coefficients, and $m_n \ll p_n$. We discuss the choice of the random matrix $R_n$ in Section \ref{sec:2.1}, and illustrate the method in detail in Section \ref{sec:2.2}. 

{\it Priors.}  We assume that the covariates are standardized. Taking a Bayesian approach to inference, we assign priors to $\bta$ and $\sigma^2$. The vector of compressed regression parameters $\bta$ is assigned a $N_{m_n}({\bf 0}, \sigma^2 I)$ prior given $\sigma^2$, where ${\bf 0}$ is a vector of zeros, and $I$ is the identity matrix. The scale parameter $\sigma^2$ is assigned an \emph{Inv-Gamma ($a_{\sigma}, b_{\sigma}$)} prior, with $a_{\sigma}, b_{\sigma}>0$.  The Normal-Inverse Gamma (ING) prior is a common choice of prior for GLMs. In the special case of a Gaussian likelihood, this prior is conjugate, and the posterior and predictive distributions are available in analytic forms.

\subsection{Choice of the projection matrix}\label{sec:2.1} 
The projection matrix $R_n$ embeds $X_n$ to a lower dimensional subspace.  If $p_n$ is not large, the best linear embedding can be estimated using a singular value decomposition (SVD) of $X_n$. 
The projection of $X_n$ to the space spanned by the singular vectors associated with the first $m_n$ singular values is guaranteed to be closer to $X_n$ in an appropriate sense than any other $m_n$-dimensional matrix.  However, if $p_n$ is very large with $p_n \gg n$, then it is problematic to estimate the projection, both computational and statistically, and random projection (RP) provides a practical alternative. If an appropriate RP matrix is chosen, due to Johnson-Lindenstrauss (JL) type embedding results, distances between sample points are maintained (see \cite{dasgupta2003}, \cite{Achlioptas_2003}). 

Our focus is on modifying current approaches by constructing RP matrices that incorporate sparsity in a way that the predictors $x_j$ having relatively weak marginal relationships with $y$ are less likely to be included in the matrix. In particular, the TArgeted Random Projection (TARP) matrices are constructed as follows: 
\begin{eqnarray}
\boldsymbol{ \gamma}&=&\left( \gamma_1, \gamma_2, \ldots, \gamma_{p_n} \right)^{\prime}\quad \mbox{and}\quad \gamma_j
\stackrel{i.i.d.}{\sim} ~\mbox{Bernoulli} \left(q_j \right)  \mbox{~~~where }\notag\\
&&\hskip17pt q_j \propto |r_{x_j,y}|^\delta~~~~ \mbox{for some constant~~ $\delta>0$,}  \label{eq_rp2} \\
R_{\overline{\bgm}} &=& O_{m_n\times (p_n - p_{\bgm})} \quad \mbox{and }~~~\Rg = R_n^{*}, \notag
 \end{eqnarray}
where $r_{x_j,y}$ is the Pearson's correlation coefficient of ${\bf x}_j$ and $\yn$, $q_j$ is the inclusion probability of $x_j$, 
$\Rg$ and $R_{\overline{\bgm}}$ are the sub-matrices of $R_n$ with columns corresponding to non-zero and zero values of $\bgm$, respectively, and $R_n^{*}$ is the $m_n\times p_{\bgm}$ projection matrix where $p_{\bgm}=\sum_j \gamma_j$.

We prioritize predictors based on their marginal utilities, ${\bf q}=(q_1,q_2,\ldots,q_{p_n})^{\prime}$, and consider a random subset of the predictors with inclusion probabilities ${\bf q}$.  This can be viewed as a randomized version of independent screening. The selected subset is further projected to a lower dimensional sub-space using $R_n^{*}$.  There are many possible choices of $R_n^{*}$ which can successfully reduce the dimension of the selected variables while having minimal impact on prediction accuracy.  Two predominant classes of such projection matrices are based on partial SVD and random projections facilitating JL type embedding. We consider both these choices, as described below.

\vskip5pt
\noindent{\bf Random projection.}
Each element, $R_{k,j}^{*}$, of $R_n^{*}$  is sampled independently from a three point distribution as 
\begin{eqnarray}
  \quad R_{k,j}^{*}&=& \left\{\begin{array}{ll}
               \pm 1/\sqrt{2\psi} \quad &\mbox{with probability}~~ \psi,\\
                    ~0 \quad &\mbox{with probability}~~ 1-2\psi, 
                    \end{array}\right. \label{eq_rp18}
\end{eqnarray}
where $\psi \in (0,0.5)$ is a constant.

Projection matrices of this form are widely used due to their inter point distance preservation property. Incorporating zero values facilitates data reduction and improves computational efficiency.   We refer to the method that
generates $R_n^*$ in  (\ref{eq_rp2}) from (\ref{eq_rp18}) as RIS-RP (Randomized Independent Screening-Random Projection).   
\begin{remark}\label{rm3}
The choice of projection matrix in (\ref{eq_rp18}) can be replaced by a wide variety of matrices having i.i.d. components with mean zero and finite fourth moments. One of the sparsest choices is of the form $R_{k,j}^{*}=\pm n^{\kappa/2}/\sqrt{m_n}$ with probability $1/2n^{\kappa}$, $0$ with probability $(1-1/n^{\kappa})$,
where $m_n\sim n^\kappa$ (see \cite{LHC_2006}). 
This choice of projection matrix reduces the data to a great extent, and is useful in compressing extremely large dimensional data. Our theoretical results would hold also if we consider a random matrix $R_n^{*}$ as above. 
\end{remark}

\vskip5pt
\noindent{\bf Principal component projection.}
Another alternative is to use the matrix of principal component scores. Let $X_{\bgm}$ be the sub-matrix of $X_n$ with columns corresponding to non-zero values of $\bgm$. Consider the partial spectral decomposition $X_{\bgm}^{\prime} X_{\bgm} = V^{\prime}_{\bgm, m_n} D_{\bgm, m_n} V_{\bgm, m_n}$. The projection matrix $R_n^{*}$ we consider is
 \vspace{-.12 in}
\begin{eqnarray}
R_n^{*} = V_{\bgm,m_n}^{\prime}. \label{eq_rp16}
\end{eqnarray}
The method corresponding to this choice of projection matrix combines a randomized version of independence screening with principal component regression (PCR). Therefore, we refer to this method as RIS-PCR (Randomized Independence Screening-PCR).

\vskip10pt
The performance of TARP depends on tuning parameters $m_n$, $\delta$ and $\psi$.  In addition, for any given choice of tuning parameters, different realizations of the random projection matrix will vary, leading to some corresponding variation in the results.  To limit dependence of the results on the choice of tuning parameters and random variation in the projection, we take the approach of generating multiple realizations of the matrix for different choices of tuning parameters and aggregating these results together.  Potentially, one could estimate weights for aggregation using Bayesian methods (see \cite{hoeting1999}) or other ensemble learning approaches, but we focus on simple averaging due to its computational and conceptual simplicity.  

\subsection{Posteriors and Predictive Distribution}\label{sec:2.2}
We illustrate the proposed method in normal linear regression for simplicity, although it is applicable to more general settings.  In the compressed regression framework, we replace the normal linear model $ y_i ={\bf x}_i^{\prime}\bbta+e_i$ by $ y_i =\left(R_n {\bf x}_i\right)^{\prime}\bta+e_i$, where $e_i\sim N(0, \sigma^2)$. Given the NIG prior structure stated above, the posterior distribution of $\bta$
follows a scaled $m_n$-variate $t$ distribution with degrees of freedom (d.f.) $n+2 a_\sigma$, location vector $\boldsymbol{\mu}_t$, and scale matrix $\Sigma_t$, where $\boldsymbol{\mu}_t=Z (X_n R_n^{\prime} )^{\prime} \yn$, $\Sigma_t=(\yn^{\prime}\yn - \boldsymbol{\mu}_t^{\prime} Z^{-1} \boldsymbol{\mu}_t+2 b_{\sigma})Z/(n+2a_\sigma)$ and $Z= (R_n X_n^{\prime} X_n R_n^{\prime} + I )^{-1}$. Moreover, the posterior distribution of $\sigma^2$, given $\mathcal{D}^n$ and $R_n$,
is inverse gamma with parameters $a_\sigma+ n/2$ and $(\yn^{\prime}\yn - \boldsymbol{\mu}_t^{\prime} \Sigma_t^{-1} \boldsymbol{\mu}_t+b_{\sigma})/2$.

Consider the problem of point prediction of $y$ when $n_{new}$ new data points on the predictors are obtained, given the dataset $\mathcal{D}^n$. The predicted values of $y$, say ${\bf y}_{new}$, can be obtained using the Bayes estimator of $\bta$ under squared error loss as follows

\vskip5pt
$ \hat{\bf y}_{new} = X_{new} R_n^{\prime} \hat{\bta}_{Bayes} ~~\mbox{where}~~ \hat{\bta}_{Bayes} =\left(Z_n^{\prime}Z_n+I\right)^{-1}Z_n^{\prime} \yn~~ \mbox{and}~~ Z_n=X_nR_n.$

\noindent Here $X_{new}$ is the new design matrix.
Moreover, the posterior predictive distribution of ${\bf y}_{new}$ is
a $n_{new}$-variate $t$ distribution with degrees of freedom $n+2 a_{\sigma}$, location vector $\hat{{\bf y}}_{new}$ and scale parameter $(\yn^{\prime}\yn - \boldsymbol{\mu}_t^{\prime} Z^{-1} \boldsymbol{\mu}_t+2 b_{\sigma})(I+X_{new}ZX_{new}^{\prime})/(n+2a_\sigma)$.

When the distribution $\yn$ is non-normal, analytical expressions of the posteriors of $\bta$, $\sigma^2$ and predictive distribution of ${\bf y}_{new}$ are not available.  In such cases, it is common to rely on a Laplace approximation (see \cite{tierney1986}) or sampling algorithm, such as MCMC. In the compressed regression framework, as the $p_n$-dimensional variables are projected to a much lower-dimensional $m_n$-hyperplane, with $m_n \ll p_n$, we are no longer in a high-dimensional setting.  Hence, MCMC is computationally tractable.

\subsection{Tuning Parameter Choice}\label{sec:3.1}
Next, we describe the choices of tuning parameters that are involved in TARP.

\noindent {\bf Choice of $m_n$.} The parameter $m_n$ determines the number of linear combinations of predictors we consider. Instead of choosing a fixed value of $m_n$, we consider a range of values.
In particular, we suggest choosing values in the range $(2\log p_n, \min\{3n/4,p_n\})$, consistent with our 
theoretical results in Section \ref{sec:3} requiring that $m_n < n$ and with our experiments assessing performance for different choices in a variety of scenarios.

\noindent {\bf Choice of $\delta$ and ${\bf q}$.} The parameter $\delta$ plays an important role in screening.  Higher values of $\delta$ lead to fewer variables selected in the screening step. 
If $p_n \gg n$, one would like to select a small proportion of predictors, and if $p_n\sim n$ then selection of a moderate proportion is desirable. We recommend  $\delta=\max\{0,(1+\log(p_n/n))/2\}$ as a default. This function selects all the predictors if $p_n\ll n$, and becomes more restrictive as $p_n$ becomes larger than $n$. The selection probabilities in the RIS stage are then $q_j=|r_{x_j,y}|^{\delta}/\max_j |r_{x_j,y}|^{\delta}$, $j=1,2,\ldots,p_n$. Hence, the variable with highest marginal correlation is definitely included in the model.

\noindent{\bf Choice of $\psi$.} The value of $\psi$ controls sparsity in the random matrix, and it is necessary to let $\psi \in (0,0.5).$  \cite{Achlioptas_2003} suggest choosing $\psi=1/6$ as a default value.
To avoid sensitivity of the results for a particular choice of $\psi$, we choose $\psi\in (0.1,0.4)$ avoiding the very sparse and dense cases.

\subsection{Computational Algorithm and Complexity}\label{sec:3.2}
We now illustrate the time complexity of TARP along with the algorithm for computation in normal linear models.

\noindent{\bf RIS-RP.} For a specific choice of $(m_n, \delta, \psi)$, calculation of $\hat{\bf y}_{new}$ using RIS-RP involves the following steps:
\begin{enumerate}[1:]
 \item Calculate $r_{x_j,y}$ for $j=1,\ldots,p_n$.
 \item Generate $\gamma_j\sim \mbox{Bernoulli}(q_j)$ where $q_j=|r_{x_j,y}|^{\delta}/\max\{|r_{x_j,y}|^{\delta}\},$ $j=1,\ldots,p_n$. \linebreak
         IF $\gamma_j=1$, generate $R_n$ with $R_{i,j}$ as in (\ref{eq_rp18}). ~ 
ELSE set $R_{i,j}=0$.
\item Post-multiply $R_n$ with $X_n$. Set $Z_n=X_nR_n$.
\item For a given $X_{new}$, compute $Z_{new}=X_{new}R_n$ and $\hat{\bf y}_{new}=Z_{new} \hat{\bta}$.
\end{enumerate}

 
 \noindent The complexity of steps 1, 2-3, 4 and 5 are $O(p_n)$, $O(n\pg m_n)$, $O(n m_n^2)$ and $O\left(n_{new} \pg m_n \right)$, respectively, where $\pg =\sum \gamma_j$. Thus, if $n_{new}\leq n$, the total complexity for a single choice of $(m_n, \delta, \psi)$ is $O(p_n)+2O(n m_n\pg )+O( n m_n^2) $ without using parallelization.

\vskip5pt
\noindent{\bf RIS-PCR.} 
RIS-PCR differs from RIS-RP in step 2 of the algorithm. After generation of $\bgm$ RIS-PCR requires SVD of $X_{\bgm}$ involving complexity $O\left(n\pg \min\{n,\pg\}\right)$. Thus, total complexity of RIS-PCR for a single choice of $(m_n, \delta, \psi)$ can similarly be derived. Therefore, the two methods have comparable time complexity unless either $n$ or $\pg$ is much larger than $m_n$. Although theoretically we do not impose any restriction on $\pg$, in practice when $p_n=\exp\{o(n)\}$ and $\delta\geq 2$, $\pg$ is usually of order $n$. 

\vskip5pt
\noindent{\bf Increment of complexity due to aggregation.}
Suppose $N$ different choices of $(m_n, \psi, R_n)$ are considered. Each choice yields a model $\mathcal{M}_l : y\sim f\left(y | {\bf x}, m_{n,l}, \psi_l, R_{n,l} \right)$ along with a corresponding estimate of ${\bf y}_{new}$ (say $\hat{\bf y}_{new,l}$), where $l\in\{1,2, \ldots, N\}$. The proposed estimate is the simple average of these $N$ estimates of ${\bf y}_{new}$.

Step 1 in the TARP algorithm is not repeated over the aggregration replicates, while the remaining steps are repeated $N$ times.  In addition, the first step of screening and aggregration are embarassingly parallelizable.  Hence, given $k$ processors, if $n_{new}\leq n$, the total complexity is $O(p_n/k)+2O(Nn m_n\pg/k )+ O(N n m_n^2/k) $ for RIS-RP and $O(p_n/k)+O( N n \pg \min\{n, \pg\}/k) +O(Nn m_n\pg/k )+ O( N n m_n^2/k )$ for RIS-PCR.

\section{Theory on Predictive Accuracy} \label{sec:3}

We study the asymptotic performance of the predictive distribution produced by TARP for a single random projection matrix without considering the aggregation step.  We focus on \emph{weakly sparse} and \emph{dense} cases where the absolute sum of the true regression coefficients is bounded. This condition also includes \emph{strong sparsity} where only a few covariates have non-zero coefficients. 

The projection matrix in TARP depends on the random variable $\bgm$, and therefore is denoted by $\Rg$. 
We denote a particular realization of the response variable as $y$, and a particular realization of the variables $(x_1, x_2, \ldots, x_{p_n})^{\prime}$ as ${\bf x}$.  
Let $f_0$ be the true density of $y$ given the predictors, and $f(y|{\bf x}, \bgm ,\Rg, \bta)$ be the conditional density of $y$ given the model induced by $\left(\bgm,\Rg\right)$, the corresponding vector of regression coefficients $\bta$ and ${\bf x}$. We follow \cite{Jiang_2007} in showing that the predictive density under our procedure is close to the true predictive density in an appropriate sense.

We assume that each covariate $x_j$ is standardized so that $|x_j|<M$, for $j=1,2,\ldots,p_n$, with $M$ a constant. We also assume that the scale parameter $\sigma^2$ in (\ref{eq_rp1}) is known.  We require the following two assumptions on the design matrix.

\vskip5pt
{\bf Assumption (A1)}  Let $r_{x_j,y}$ denote the correlation coefficient of the observed values of $x_j$ and $y$, $1\leq j\leq p_n$. Then for each data point $(y,{\bf x})$ and constant $\delta$ in (\ref{eq_rp2}), there exists a positive constant $\alpha_{\delta}$ such that
\vspace{-.1 in}
\[  \lim_{p_n\rightarrow\infty} \frac{1}{p_n} \sum_{j=1}^{p_n} x_{j}^2 |r_{x_j,y}|^{\delta} \rightarrow {\alpha}_{\delta}. \]
\vspace{-.1 in}

{\bf Assumption (A2)} 
Let $q(\bgm)= \prod_{i=1}^{n} q_j^{\gamma_j} (1-q_j)^{(1-\gamma_j)}$, with $q_j$ defined in (\ref{eq_rp2}), denote the probability of obtaining a particular $\bgm=(\gamma_1,\ldots,\gamma_{p_n})'$ in the random screening step.  Let $\Gamma_l \subset \{0,1\}^{p_n}$ denote the set of $\bgm$ vectors such that $p_{\bgm} = \sum_{j=1}^{p_n}\gamma_j = l$, and let $\mathcal{M}_l \subset \Gamma_l$ denote the first $p_n^{k_n}$ elements of $\Gamma_l$ ordered in their $q(\bgm)$ values.  Let $\mathcal{A}_n$ denote the event that $\gamma \in \mathcal{M}_l$ for some $l$.  Then, 
$P(\mathcal{A}_n^c)=P\big( \{\bgm: \bgm \notin \cup_l \mathcal{M}_l\} \big)
\le \exp(-n\varepsilon_n^2/4),$
for some increasing sequence of integers $\{k_n\}$ and sequence $\{\varepsilon_n\}$ satisfying $0<\varepsilon_n^2<1$ and $n\varepsilon_n^2\rightarrow\infty$.

\begin{remark}\label{rm2}
As the probability of selection in the random screening step depends on the empirical correlation between the predictor and the response, assumption (A2) is on the data generating process.  If $l \le k_n$ or $l \ge p_n-k_n$, then all models of dimension $l$ belong to $\mathcal{M}_l$, and none of the corresponding $\gamma$ vectors belong to $\mathcal{A}_n^c$.  If $k_n < l < p_n-k_n$, there are more than $p_n^{k_n}$ models of dimension $l$, but as the models are ordered in terms of their $q({\bgm})$ values, the models not falling in 
$\mathcal{M}_l$ should have extremely small probabilities of selection, hence satisfying (A2).  Violations of (A2) would imply that large numbers of predictors have empirical correlations that are not close to zero.
\end{remark}


{\it Measure of closeness:} Let $\nu_{\bf x}(d{\bf x})$ be the probability measure for ${\bf x}$, and $\nu_y(dy)$ be the dominating measure for conditional densities $f$ and $f_0$. The dominating measure of $(y,{\bf x})$ is taken to be the product of $\nu_y(dy) \nu_{\bf x}(d{\bf x})$.

The Hellinger distance between $f$ and $f_0$ is given by

  \vskip5pt
\qquad \qquad \quad \qquad \qquad $d(f,f_0)=\sqrt{\displaystyle\int \left(\sqrt{f}-\sqrt{f_0} \right)^2 \nu_{\bf x}(d{\bf x}) \nu_y(dy)}.$

 \vskip5pt
The Kullback-Leibler divergence between $f$ and $f_0$ is given by 

  \vskip5pt
\qquad \qquad \quad \qquad \qquad $ d_0(f,f_0)=\displaystyle\int f_0 \ln \left(\frac{f_0}{f} \right) \nu_{\bf x}(d{\bf x})\nu_y(dy).$

\vskip5pt
Define  $d_t(f,f_0)=t^{-1}\Big( \displaystyle\int f_0 \displaystyle\left(\frac{f_0}{f} \right)^t\nu_{\bf x}(d{\bf x})\nu_y(dy) -1\Big),$ for any $t>0.$

\vskip5pt
Consider the following two facts: (i) $d(f,f_0) \leq \left( d_0(f,f_0) \right)^{1/2}$, and (ii) $d_t(f,f_0)$ decreases to $d_0(f,f_0)$ as $t$ decreases to $0$  (see \cite{Jiang_2007}).

Let $\mathcal{P}_n$ be a sequence of sets of probability densities, and $\varepsilon_n$ be a sequence of positive numbers.
Let $ N(\varepsilon_n,\mathcal{P}_n)$ be the $\varepsilon_n$-covering number, i.e., the minimal number of Hellinger balls of radius $\varepsilon_n$ needed to cover $\mathcal{P}_n$. 

\paragraph{RIS-RP.} The result showing asymptotic accuracy in approximating the predictive density using RIS-RP is stated below.

\begin{theorem}\label{thm:1}
 Let $\bta\sim N({\bf 0},\sigma_{\theta}^2 I)$, and $f(y|{\bf x}, \bgm ,\Rg, \bta)$ be the conditional density of $y$ given the model induced by $\left(\bgm,\Rg\right)$, where $\Rg$ is as in (\ref{eq_rp2}) and (\ref{eq_rp18}). Let $\bbta_{0}$ be the true regression parameter with $\sum_j |\beta_{0,j} |<K$ for some constant $K$, and assumptions (A1)-(A2) hold. Consider the sequence $\{ \varepsilon_n\}$ as in assumption (A2) satisfying $0< \varepsilon_n^2<1$ and $n\varepsilon_n^2 \rightarrow \infty$, and assume that the following statements hold for sufficiently large $n$:\\
(i) $ m_n |\log \varepsilon_n^2| < n \varepsilon_n^2/4$, \\
(ii) $ k_n \log p_n < n \varepsilon_n^2/4$, and \\
(iii) $m_n \log \left(1 + D\left(\sigma_{\theta} \sqrt{6 n \varepsilon_n^2 p_n m_n} \right) \right) < n \varepsilon_n^2/4$, where

\qquad $D(h^{*})=h^{*} \sup_{h\leq h^{*}} |a^{\prime} (h)| \sup_{h\leq h^{*}} |a^{\prime} (h)/b^{\prime}(h)|$, $b(\cdot)$ as in (\ref{eq_rp1}). Then, $$P_{f_0}\left[\pi\left\{d(f,f_0)> 4\varepsilon_n| \mathcal{D}^n \right\} > 2 e^{-n\varepsilon_n^2/4 } \right] \leq 2 e^{-n\varepsilon_n^2/5 },  $$
where $\pi\{ \cdot | \mathcal{D}^n\}$ is the posterior measure.
\end{theorem}
\noindent The proof of Theorem \ref{thm:1} is given in the Appendix.

\begin{remark}\label{rm4}
 If we consider the sparse choice of $\Rg$, as described in Remark \ref{rm3}, the same line of proof (as that of Theorem \ref{thm:1}) would go though. For both the choices of $\Rg$, each component of the projection matrix has expectation zero and finite fourth moment.
For the random matrix in (\ref{eq_rp18}), the probability of choosing a non-zero element, $P(R_{i,j}\neq 0)$, is fixed, while the probability is decaying with $n$ for the sparse choice of $\Rg$. However, the rate of decay is such that distances are preserved between the sample-points, a critical property for proving consistency.
\end{remark}

\paragraph{RIS-PCR.} Asymptotic guarantees on predictive approximation for RIS-PCR requires 
 an additional assumption.

\vskip5pt
\noindent{\bf Assumption (A3)} Let $X_{\bgm}$ be the sub-matrix of $X_n$
with columns corresponding to non-zero values of $\bgm$, and ${\bf x}_{\bgm}$ be a row of $X_{\bgm}$. Let $V_{\bgm}$ be the $m_n \times \pg$ matrix of $m_n$ eigenvectors corresponding to the first $m_n$ eigenvalues of $X_{\bgm}^{\prime} X_{\bgm}$.
Then, for each $\bgm$ and data point ${\bf x}_{\bgm}$,
$$ \frac{\| V_{\bgm} {\bf x}_{\bgm} \|^2}{\| {\bf x}_{\bgm} \|^2} \geq \alpha_n,$$
where $\alpha_n \sim (n\varepsilon_n^2)^{-1}$, where the sequence $\{ \varepsilon_n^2\}$ is as in assumption (A2).

\begin{remark}\label{rm5}
If the matrix $X_{\bgm}^{\prime} X_{\bgm}$ has rank less than $m_n$, then $\alpha_n=1$ by Perseval's identity. Consider the situation where rank of the gram matrix, say $r_n (\leq n)$, is bigger than $m_n$. Then the row space of $X_{\bgm}$, or that of $X_{\bgm}^{\prime} X_{\bgm}$, is spanned by a set of $r_n$ basis vectors ${\bf v}_1, {\bf v}_2, \ldots, {\bf v}_{r_n}$. Therefore, any data point ${\bf x}$ can be written as a linear combination of these $r_n$ vectors as 
${\bf x}=a_1 {\bf v}_1+ a_2 {\bf v}_2+ \cdots + a_{r_n} {\bf v}_{r_n}$,
where $a_1, a_2, \ldots, a_{r_n}$ are constants not all equal to zero. As the vectors ${\bf v}_j$ are orthonormal, ${\bf v}_j^{\prime} {\bf x}=a_j$ for all $j=1,2, \ldots, r_n$, which in turn implies that ${\bf x}^{\prime} {\bf x}=\sum_{j=1}^{r_n} a_j^2$. Also, note that the first $m_n$ among these $r_n$ vectors constitute $V_{\bgm}^{\prime}$, which implies $\|V_{\bgm}^{\prime} {\bf x} \|^2=\sum_{j=1}^{m_n} a_j^2$. Thus $\| V_{\bgm} {\bf x} \|^2/\| {\bf x} \|^2= \sum_{j=1}^{m_n} a_j^2 / \sum_{j=1}^{r_n} a_j^2$, and magnitude of the ratio depends on the part of ${\bf x}$ explained by the last few principal component directions. The lower bound $\alpha_n \sim (n\varepsilon_n^2)^{-1}$ is weaker than many real data 
scenarios where most of the variation is explained by the first few principal components.
\end{remark}

\begin{theorem}\label{thm:2}
 Let $\bta\sim N({\bf 0},\sigma_{\theta}^2 I)$, and $f(y|{\bf x}, \bgm ,\Rg, \bta)$ be the conditional density of $y$ given the model induced by $\left(\bgm,\Rg\right)$, where $\Rg$ is as in (\ref{eq_rp2}) and (\ref{eq_rp16}). Let $\bbta_{0}$ be the true regression parameter with $\sum_j |\beta_{0,j} |<K$ for some constant $K$, and assumptions (A1)-(A3) hold. Assume that the conditions (i)-(iii) of Theorem \ref{thm:1} hold for the sequence $\{ \varepsilon_n\}$ as in assumption (A2) satisfying $0< \varepsilon_n^2<1$ and $n\varepsilon_n^2 \rightarrow \infty$.
%
Then, $$P_{f_0}\left[\pi\left\{d(f,f_0)> 4\varepsilon_n| \mathcal{D}^n \right\} > 2 e^{-n\varepsilon_n^2/4 } \right] \leq 2 e^{-n\varepsilon_n^2/5 },  $$
where $\pi\{ \cdot | \mathcal{D}^n\}$ is the posterior measure.
\end{theorem}
\noindent The proof of Theorem \ref{thm:2} is given in the Appendix.

\begin{remark}\label{rm6} The conditions (i)-(iii) in Theorems \ref{thm:1} and \ref{thm:2} are related to the sizes of $p_n$, $m_n$ and $k_n$ in comparison with $n\varepsilon_n^2$.
A sufficient condition for (i) is $m_n \log n < n\varepsilon_n^2/4$, providing an upper bound on the dimension of the subspace $m_n$. Condition (ii) restricts the permissible number of regressors $p_n$, and the number of possible models of each dimension $k_n$. If there is a strict ordering in the marginal correlation coefficients $|r_{x_j,y}|$, so that $k_n\leq \kappa$ for some large number $\kappa$ (see assumption (A2)), then the condition reduces to $\log p_n < n\varepsilon_n^2/4$.  To illustrate that condition (iii) tends to be weak, consider distributions of $y$ corresponding to Bernoulli, Poisson and normal. For these cases, the quantity $D(h^{*})$ is at most order $O(h^{*})$. Therefore, condition (iii) does not impose much additional restriction over (i)-(ii), except $m_n \log p_n < n\varepsilon_n^2/4$, inducing a stronger upper-bound to $m_n$.  
\end{remark}

\section{Simulation Study} \label{sec:4}
In this section, we consider different simulation schemes (\emph{Scheme I} -- \emph{IV})  to compare TARP with a variety of methods. We mainly focus on high-dimensional and weakly sparse regression problems with a variety of correlation structures in the predictors.  The sample size is taken to be $200$, while $p_n$ varies.  Additional results for different choices of $n$ are provided in the Supplement. 

\vskip5pt
\noindent{\it Competitors.} 
We compare with: SCAD screened by iterative SIS (ISIS), \emph{ISIS-SCAD}; minimax concave penalty (MCP) method screened by ISIS, \emph{ISIS-MCP}; LASSO screened by sequential strong rule (SSR, \cite{SSR}), \emph{SSR-LASSO}; ridge regression screened by SSR, \emph{SSR-Ridge}; elastic net screened by SSR,  \emph{SSR-EN}; principal component regression (\emph{PCR}); sparse PCR, \emph{SPCR} (see \cite{SPCR}); robust PCR, \emph{RPCR} (see \cite{RPCA}); and Bayesian compressed regression (\emph{BCR}).
ISIS-SCAD and ISIS-MCP are available in the `SIS' package, and LASSO, ridge and elastic net are available in the ‘biglasso’ package (\cite{biglasso}). SPCR and RPCR are performed using `PMA' and `rsvd' packages in R, respectively.  To estimate PC scores, we rely on approximate SVD using \emph{fast.svd} in the `corpcor' package. For BCR, we average over $100$ different random projections with varying $m_n$ values within the range $[2\log p_n,3n/4]$. We use the \emph{qr} function in R to apply QR factorization in place of Gram-Schmidt orthogonalization of the random matrix, which is computationally prohibitive for large $p_n$.

\vskip5pt
\noindent{\it The proposed method.} 
We select the tuning parameters of TARP as described in Section \ref{sec:3.1}. The parameter $m_n$ is chosen in the range $[2\log p_n,3n/4]$. We assign $\delta=2$ as the function $\max\{0,(1+\log(p_n/n))/2\}$ is close to 2 for all the choices of $(n,p_n)$.  Further, the hyperparameters of the inverse gamma priors are set to $0.02$ to correspond to a minimally informative prior.  

\vskip5pt
\noindent{\it Simulation Schemes.}
In the first three simulation schemes, the predictors were generated from $N({\bf 0},\Sigma)$, with different choices of $p_n$, $\Sigma$ and the regression coefficients. In \emph{Scheme IV} we consider a functional regression setup. Different methods are compared with respect to their performance in out of sample prediction. We calculate mean square prediction error (MSPE), empirical coverage probability (ECP) of a $50\%$ prediction interval (PI) and the width of the PI for each of a 100 replicate datasets in each simulation case. 

\afterpage{
\begin{figure}
\centering     
\includegraphics[scale=.25]{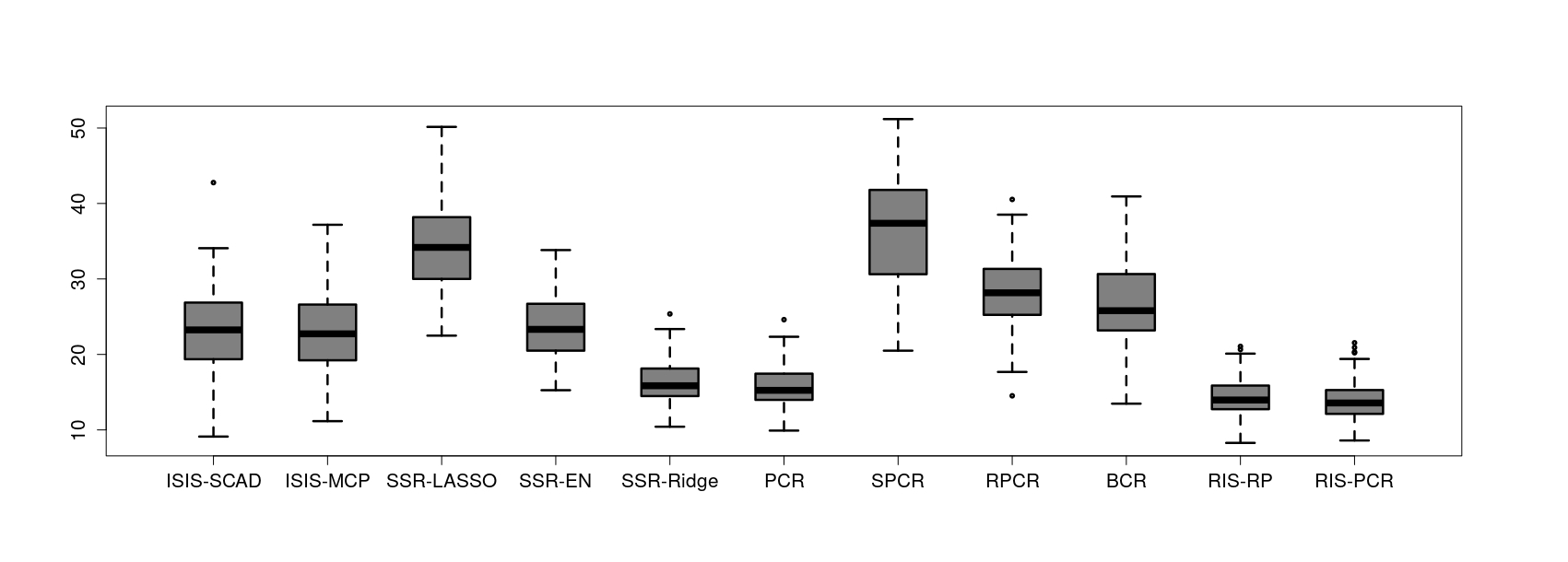}
\caption{Box-plots of MSPE for $p_n=2000$ in \emph{Scheme I}. }
\label{fig1}
\end{figure}
 }

\vskip5pt
\noindent{\it Scheme I: First order autoregressive structure.} $\Sigma_{i,j}=(0.9)^{|i-j|}$, $i,j=1,\ldots,p_n$, 
with $p_n \in \{2000,3000\}$, and $\beta_j=0$ for all but a randomly selected set of 30 predictors having
$\beta_j=1$.

 \afterpage{
\begin{figure}
\centering     
\includegraphics[scale=.25]{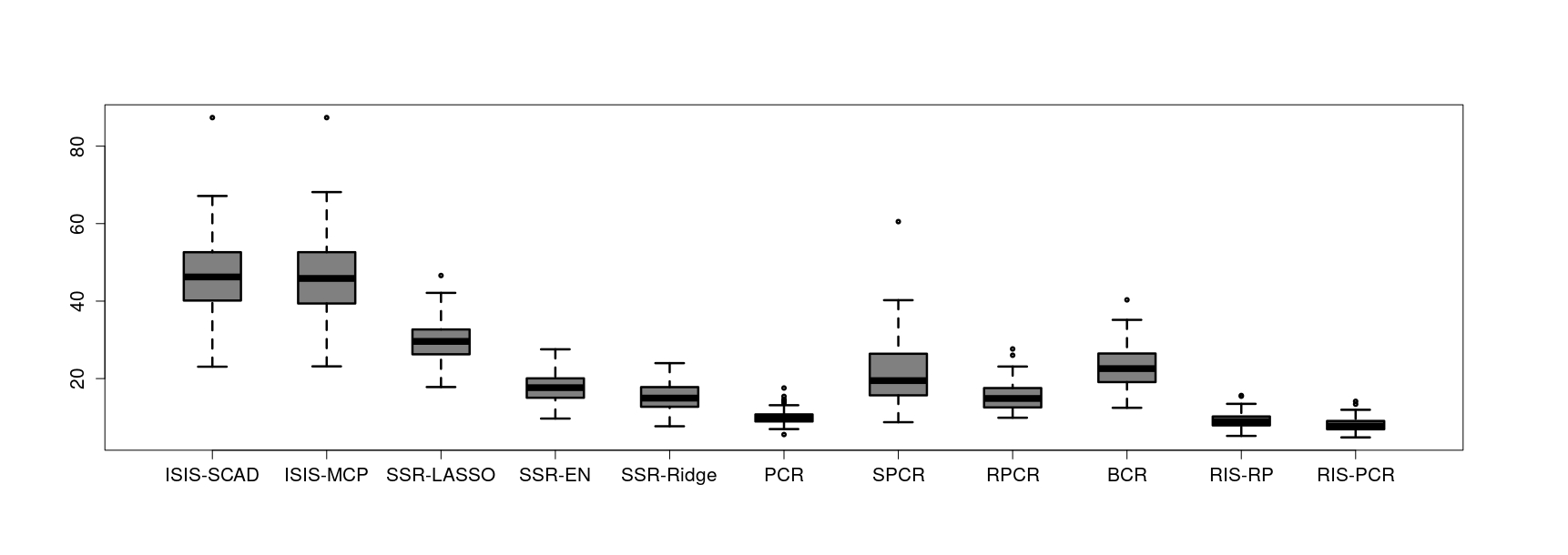}
\caption{Box-plots of MSPE for $p_n=10000$ in \emph{Scheme II}. }
\label{fig3}
\end{figure}
 }

\vskip5pt
\noindent{\it Scheme II: Block diagonal covariance structure.} We choose $(p_n/100 -2)$ blocks of $100$ predictors each, along with $200$ independent predictors, with $p_n \in \{10^4,2\times 10^4\}$. The within-block correlation is $\rho$ and the across-block correlation is zero, with $\rho=0.3$ for half of the blocks and $\rho=0.9$ for the remaining half.  There are 21 non-zero $\beta_j$s having $\beta_j=1$, with 20 of the corresponding predictors in the $\rho=0.9$ blocks and the remaining in the independent block.


\afterpage{
\begin{figure}
\centering     
\includegraphics[scale=.25]{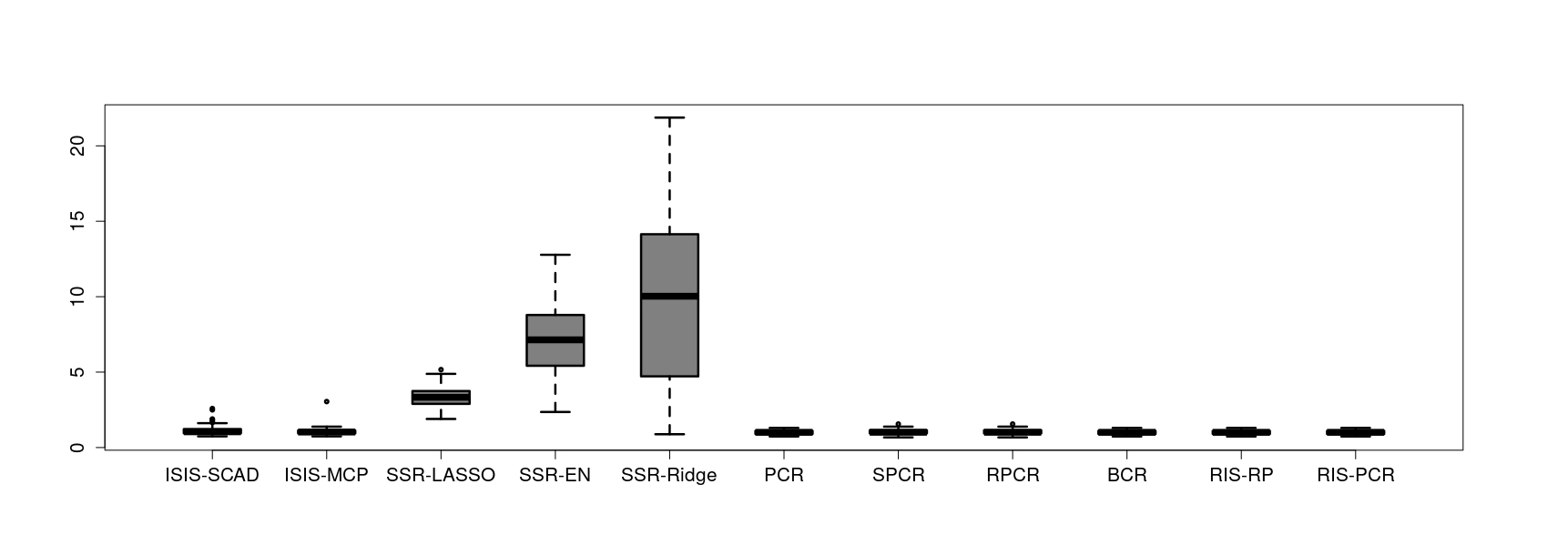}
\caption{Box-plot of MSPEs for $p_n=5000$ in \emph{Scheme III}. }
\label{fig5}
\end{figure}
 }

\vskip5pt
\noindent{\it Scheme III: Principal Component Regression.} We first choose a matrix $P$ with orthonormal columns, and a $3 \times 3$ matrix $D=\mbox{diag}(15^2,10^2,7^2)$.  We set $\Sigma=PDP^{\prime}$ and choose $\boldsymbol{\beta}=P_{\cdot,1}$, where $P_{\cdot,1}$ is the first column of $P$.  This produces an $X_n$ with three dominant principal components, with the response $\yn$ dependent on the first and
$p_n \in \{10^4,5\times 10^4\}$. 

\vskip5pt
\noindent{\it Scheme IV: Functional Regression.} Finally, we consider a functional regression setup, where the covariates are generated from Brownian bridge $B_t$ with $t\in (0,5)$ and values ranging from $(0,10)$. A set of $20$ covariates is randomly selected as active, each having regression parameters in the range $(2, 2.5)$, 
and $p_n \in \{10^4,2 \times 10^4\}$.

\vskip5pt
\noindent{\it Results.} For each of the simulation schemes, we present the results corresponding to the first choice of $p_n$, with the other results provided in the Supplement.

In \emph{Scheme I} (see Figure \ref{fig1}) SSR-Ridge, PCR, RIS-RP and RIS-PCR show competitive performance, with RIS-PCR showing the best overall result.  Performance of RIS-PCR is closely followed by RIS-RP, which in turn is followed by PCR and SSR-Ridge. Apart from these four methods, ISIS based approaches and SSR-EN exhibit reasonable performance. Although ISIS based methods have lower average MSPE, the variances of MSPEs are high indicating less stability. SPCR, RPCR, BCR and SSR-LASSO fail to perform adequately. 

\afterpage{
\begin{figure}
\centering     
\subfigure[MSPE of all the methods.]{\label{fig7}  \includegraphics[height=1.4 in, width=3 in]{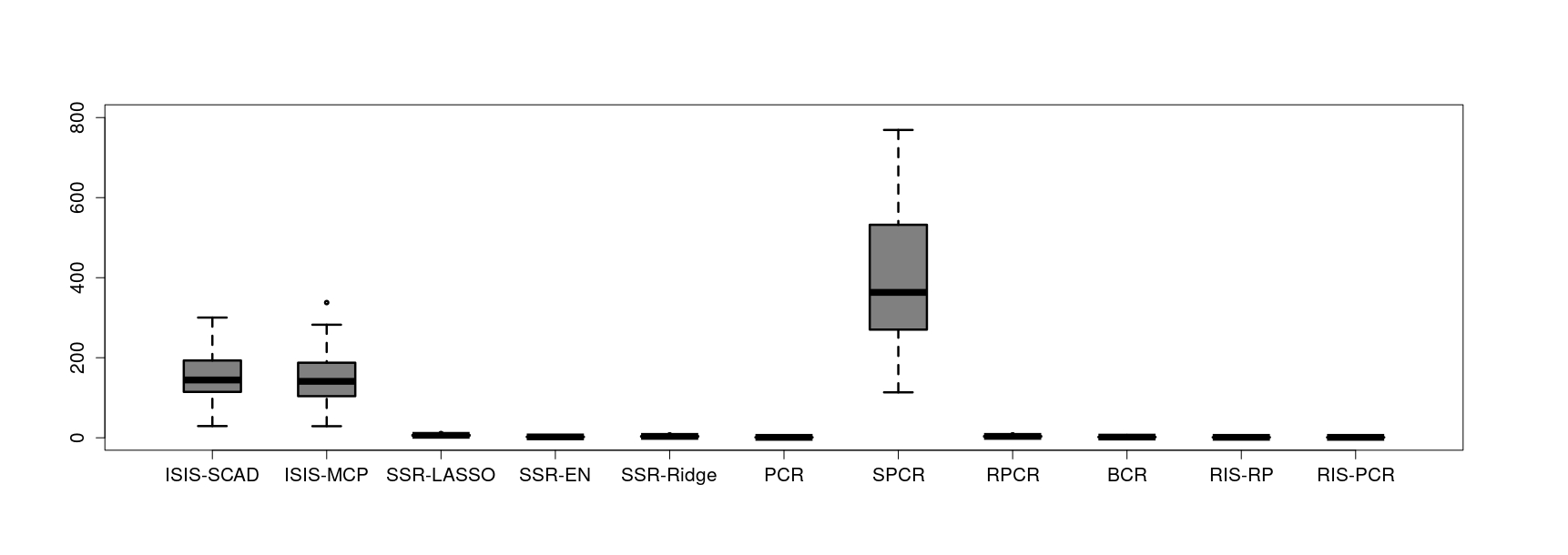}}
\subfigure[MSPE of selected methods.]{\label{fig8}  \includegraphics[height=1.4 in, width=3 in]{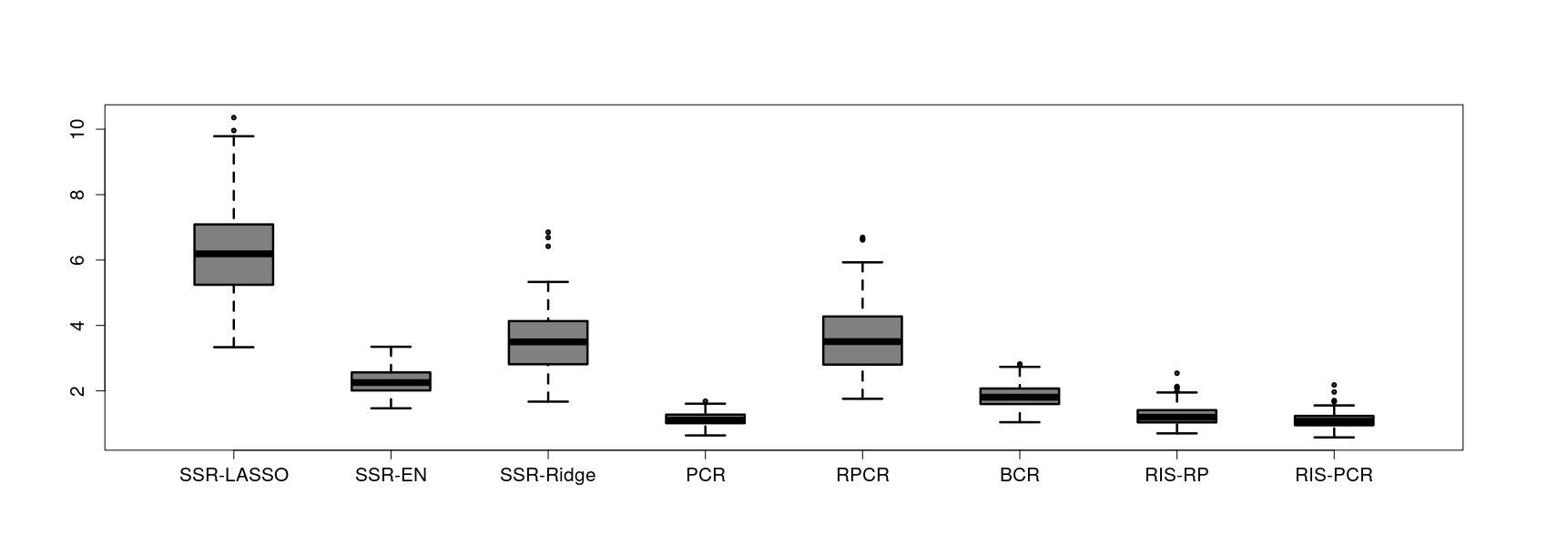}}
\caption{Box-plot of MSPEs for $p_n=10000$ in \emph{Scheme IV}. } 
\end{figure}
 }

For \emph{Scheme II}, PCR, RIS-RP and RIS-PCR again have the best performance (see Figure \ref{fig3}). RIS-PCR yields lowest overall MSPE, closely followed by RIS-RP and PCR. SSR-EN, SSR-Ridge and RPCR exhibit moderate performance. Performance of these three methods are followed by LASSO, SPCR and BCR.  Although SPCR enjoys better average MSPE, it yields higher dispersion as well. ISIS based methods fail to perform adequately.

In \emph{Scheme III} (see Figure \ref{fig5}) SSR based methods fail to exhibit competitive performance. Among the SSR based methods, performance of SSR-LASSO is better than the others. All the other methods perform well, although ISIS based methods occasionally show higher MSPEs. 

In \emph{Scheme IV} (see Figure \ref{fig7}), SPCR and ISIS based methods fail completely, making the other
methods indistinguishable in the figure.  Hence, we separately show these methods in Figure \ref{fig8}.
Among the other methods, PCR, RIS-RP and RIS-PCR have the best overall performance closely followed by BCR and then by SSR-EN. The other three methods fail to exhibit competitive performance.

\vskip10pt

We next consider the empirical coverage probabilities (ECPs) of $50\%$ prediction intervals (PI), and the width of PIs. For the Bayesian methods, the PI is obtained from the highest credible region of the predictive distribution of ${\bf y}_{new}$ given $\mathcal{D}^n, X_{new}$. The PIs for the frequentist methods can be obtained as $y_{new} \pm t_{n-p_{\bgm}-1,\frac{\alpha}{2}} \sqrt{MSPE \left(1+ {\bf x}_{\bgm,new}^{\prime} (X_{\bgm}^{\prime}X_{\bgm})^{-1}{\bf x}_{\bgm,new}\right)}$, where $t_{n-p_{\bgm}-1,\frac{\alpha}{2}}$ is the upper $\alpha/2$ point of $t$ distribution with $(n-p_{\bgm}-1)$ degrees of freedom, the suffix $\bgm$ indicates consideration of the regressors selected by the corresponding method, and $p_{\bgm}$ is the number of selected regressors. The PIs for PCR based methods can be obtained similarly, except ${\bf x}_{\bgm,new}$ and $X_{\bgm}$ are replaced by the principal component scores.
For SSR-Ridge and SSR-EN, $p_{\bgm}$ is much larger than $n$. Therefore, the PIs of these methods could not be calculated using the above formula, instead we obtain the interval as $y_{new} \pm t_{n,\frac{\alpha}{2}} \sqrt{MSPE +se(\hat{y}|\mathcal{D}^n,{\bf x}_{new})^2},$ or $y_{new} \pm t_{n,\frac{\alpha}{2}} \sqrt{2MSPE}$, whichever gives better results, where $se(\hat{y}|\mathcal{D}^n,{\bf x}_{new})^2$ is the variance of the fitted values. 
The results are summarized in Table \ref{tab1}. 
 
 \afterpage{
\begin{table*}[!h]
\renewcommand{\arraystretch}{1.5}
\begin{center}
{\scriptsize
\caption{Mean and standard deviation(sd) of empirical coverage probabilities of 50\% prediction intervals, and mean and sd of width of 50\% prediction intervals.}
\label{tab1}
\begin{tabular}{|p{1.6 cm}|p{0.7 cm} p{0.7 cm} p{0.7 cm} p{0.7 cm} p{0.7 cm} p{0.7 cm} p{0.7 cm} p{0.7 cm} p{0.7 cm} p{0.7 cm} p{0.8 cm}|}  \hline
Methods$\rightarrow$ \quad Scheme,~ $p_n$  &  ISIS-SACD &  ISIS-SACD &  SSR-LASSO &  SSR-EN &  SSR-Ridge & \qquad PCR & \qquad SPCR & \qquad RPCR & \qquad BCR & RIS-RP & RIS-PCR \\ \hline 
\multicolumn{12}{|c|}{\emph{Average and standard deviation (in braces) of empirical coverage probability}} \\ \hline
I~~   \; ~$2\times 10^3$ & 0.309 (.052) & 0.311  (.058) & 0.432 (.051) & 0.549  (.040) & 0.692 (.055) & 0.498 (.051) & 0.457 (.055) & 0.433 (.290) & 0.286 (.053)  & 0.493 (.059)  & 0.431 (.057)  \\ 

II~  \; ~$10^4$ & 0.327 (.056) & 0.324 (.055) & 0.429 (.049) & 0.607 (.062) & 0.702 (.075) & 0.455 (.308)  &0.499 (.058) & 0.445 (.050)&  0.278 (.049) &  0.502 (.056) & 0.362 (.050)  \\ 

III  \; ~$5\times 10^3$ & 0.494 (.053) & 0.494 (.053) &  0.503 (.058) &  0.657 (.031) &  0.678 (.127) & 0.503 (.058) & 0.494 (.049) & 0.494 (.049) & 0.527 (.054) & 0.494 (.055) & 0.507 (.054)  \\  

IV~  \; ~$10^4$ & 0.425 (.053) & 0.416 (.056) & 0.477 (.065) & 0.660 (.031) & 0.665 (0.102) & 0.448 (.300) &0.488 (.048) & 0.491 (.053)& 0.487 (.058) & 0.696 (.052) &  0.418 (.055)  \\
\hline
\multicolumn{12}{|c|}{\emph{Average and standard deviation (in braces) of width of the 50\% prediction interval}} \\ \hline
I~~   \; ~$2\times 10^3$ & 3.894 (.500) & 3.857 (.480) & 5.939 (.424) & 7.410 (.578) & 8.192 (.528) & 7.739 (8.309) &7.964 (.700) & 6.483 (.451) & 5.260 (.416) & 5.029 (.194) & 4.253 (.365) \\ 

II~  \; ~$10^4$ & 5.868 (.657) & 5.836 (.684) & 5.113 (.527) & 7.335 (.615) & 8.139 (.830) & 8.407 (25.713) & 6.066 (1.103) & 4.635 (.407) &  4.894 (.433) &  4.204 (.198)  & 2.744 (.232)  \\ 

III  \; ~$5\times 10^3$ & 1.424 (.174) & 1.385 (.144) & 1.362 (.072) & 1.945 (.145) & 5.755 (1.928) &  1.362 (.072)& 1.366 (.069) & 1.366 (.069) & 1.463 (.073) & 1.351 (.079) & 1.391 (.077)   \\

IV~  \; ~$10^4$ & 13.63 (2.705) & 13.34 (2.540) & 2.170 (.422) & 1.537 (.139) &  3.582 (.501) & 2.508 (3.169) & 26.512 (6.209) & 2.480 (.345) & 1.792 (.115) & 2.284 (.268) & 1.154 (.123)  \\
\hline
\end{tabular}
}
\end{center}
\end{table*}
 }

\paragraph{\it Summary.}  Considering all the simulation schemes, RIS-RP, RPCR, SPCR and SSR-LASSO have the best performance with respect to mean ECP, with RIS-RP having lowest average width among these four methods. The average width of PI is much bigger for SPCR, particularly in \emph{Scheme IV}.
ISIS based methods and RIS-PCR have relatively lower coverage probabilities in general, although among these methods, RIS-PCR has higher coverage with much lower average width than others, especially in \emph{Scheme IV}. The average ECP for PCR is satisfactory, although the corresponding widths of PIs have large variances in almost all the simulation schemes. This indicates instability in the overall performance of PCR. BCR shows under-coverage in the first two simulation schemes, but performs well with respect to both the measures in \emph{Schemes III} and \emph{IV}.  Finally the other two methods, viz., SSR-Ridge and SSR-EN have higher values of ECP, along with higher width of PIs. SSR-Ridge has highest average width of PI in \emph{Schemes I} and \emph{III}. In all the simulation schemes SSR-EN outperforms SSR-Ridge with respect to width of PI.

\vspace{-.05 in}
\subsection{Computational Time} \label{sec:4.1}
The computational time of a method may depend on the simulation scheme due to varying level of complexity  in the dataset. We only present the computational time for \emph{Scheme IV} as an example. Figure \ref{fig11} presents the time (in minutes) taken by different methods to compute $\hat{\bf y}_{new}$ using a single core, as $p_n$ grows and $n=n_{new}=100$. We run all the methods in R 3.4.1 in a 64 bit Dell-Inspiron desktop with Ubuntu 16.04 LTS operating system, 15.6 GB random access memory and Intel® Core™ i5-4460 CPU @ 3.20GHz processor.   

\afterpage{
\begin{figure}
\begin{center}     
 \includegraphics[scale=.35, angle=270]{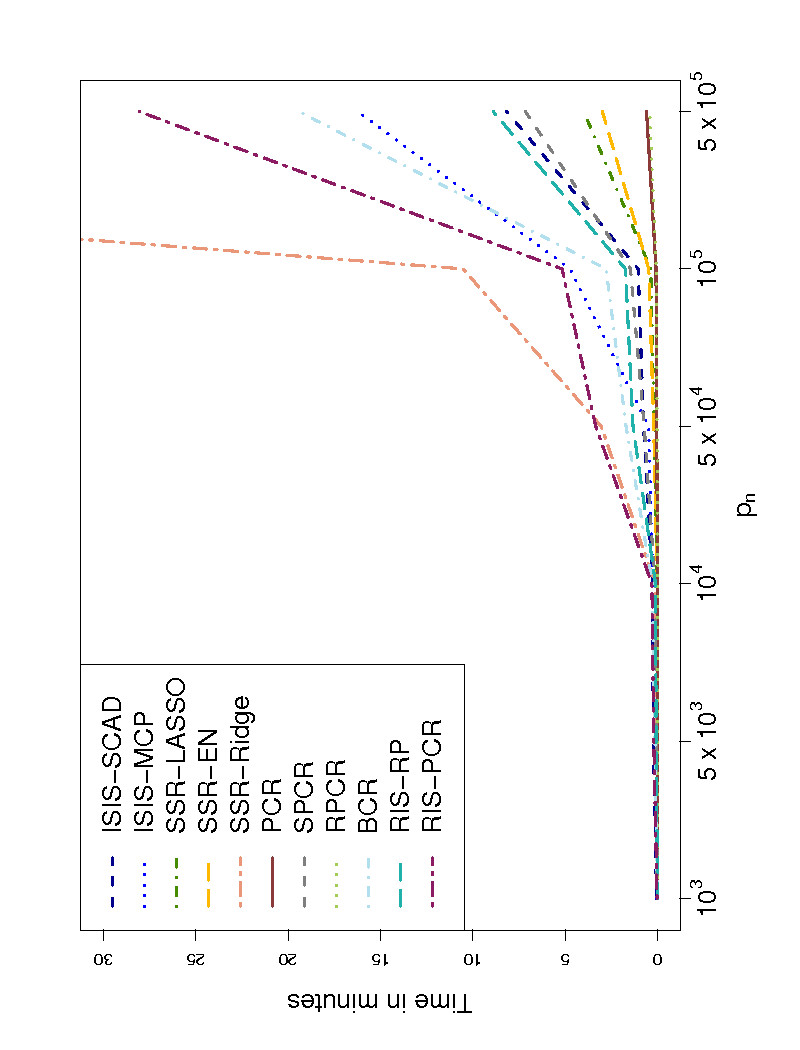}
\caption{System time required by different methods to predict $y$ as $p_n$ grows.}
\label{fig11}
\end{center}
\end{figure}}
 {\it Results.} When $p_n$ is below $10^4$ all the methods require comparable computational time. When $p_n$ is increased, SSR based methods, except SSR-Ridge, RPCR, and \emph{fast.svd} based PCR continue to require low computational time. ISIS-SCAD, SPCR and RIS-RP also have reasonable computational expense (approximately, 5 minutes for $p_n=5\times 10^5$). Computational time of BCR, ISIS-MCP and RIS-PCR tends to increase rapidly after $p_n=10^5$. Among these three methods, RIS-PCR requires highest system time (approximately 27 minutes for $p_n=5\times 10^5$). The computational time required by SSR-Ridge exceeds all other methods for $p_n>5\times 10^4$, and for $p_n=5\times 10^5$ it becomes computationally prohibitive (it takes more than $2$ hours).  

The increment of computational time of RIS-PCR is due to the computation of exact SVD of the screened design matrix $X_{\bgm}$. However, this burden would immediately be reduced if one uses some approximation of the SVD. In that case the computational time would be comparable to RIS-RP.

 \section{Real Data Analysis} \label{sec:5}
 In this section, we study the performance of TARP using three real datasets, viz., \emph{Golub} dataset, \emph{GTEx} dataset and \emph{Eye} dataset. The \emph{Golub} dataset is available at GitHub \url{https://github.com/ramhiser/datamicroarray/wiki}, the \emph{GTEx} dataset is available at GTEx portal \url{https://www.gtexportal.org/home/} and the \emph{eye data} is available in the \emph{flare} package of R.  In each case, we assess out-of-sample predictive performance averaging over multiple training-test splits of the data.
 
 \vskip5pt
 \noindent{\bf Golub data.} The Golub data consist of 47 patients with acute lymphoblastic leukemia (\emph{ALL}) and 25 patients with acute myeloid leukemia (\emph{AML}). Each of the 72 ($=n$) patients had bone marrow samples obtained at the time of diagnosis (see \cite{Golub}). Expression levels of 7129 ($=p_n$)  genes have been measured for each patient. We consider a training set of size $60$ with $20$ \emph{AML} patients, and $40$ \emph{ALL} patients. The test set consists of the remaining $12$ samples.
 
 \vskip5pt
 \noindent{\bf GTEx Data.}
 To understand the functional consequences of genetic variation, \cite{gtex} presented an analysis of RNA sequencing data from 1641 samples across 43 tissues from 175 individuals, generated as part of the pilot phase of the Genotype-Tissue Expression (GTEx) project.
 We selected RNA-seq data on two normal tissues, viz., Artery-Aorta and Artery-Tibial. The dataset contains RNA-seq expressions on 36115 $(=p_n)$ genes and 556 $(=n)$ samples, among which 224 are from Artery-Aorta, and 332 are from Artery-Tibial. A training set of $100$ samples from each of the tissue types is considered, and the remaining $446$ samples are used as test set.
 
 \vskip5pt
 \noindent{\bf Eye Data.} The Eye dataset consists of gene expressions for $200$ ($=p_n$) gene probes from the microarray experiments of mammalian-eye tissue samples of 120 ($=n$) rats (see \cite{eyedata}). The response variable is the expression level of the TRIM32 gene. 
 We consider $100$ sample points as the training set, and the remaining $20$ samples as the test set.

Golub and GTEx datasets have nominal response, and therefore the methods are evaluated by the misclassification rate (in \%) and the area under the receiver operating characteristic (ROC) curve. Table \ref{tab2} provides the average and standard deviation (sd) of percentages of misclassification, and those for the area under the ROC curve over 100 random subsets of the same size chosen from the dataset for the competing methods. We further compare the predictive performance of the methods in terms of mean squared difference of predictive and empirical probabilities for these two datasets. Most methods (except SSR based methods and SPCR) exhibit similar performance in this aspect. We provide the details of the predictive calibration in the Supplement. 

The eye dataset has continuous response, and therefore we evaluate the methods by MSPE and empirical coverage probabilities (ECP) of $50\%$ prediction intervals (PI) as in Section \ref{sec:4}.  As variation in the expression levels of the TRIM32 gene is very small (the range is 1.37) we multiply the MSPEs of different methods by $10$ to increase the variability. Table \ref{tab2} provides the mean and sd of MSPEs, ECPs of $50\%$ PIs, and widths of the PIs over 100 different training and test sets selected from the dataset, for the competing methods.  
 \afterpage{
\begin{table*}[!h]
\renewcommand{\arraystretch}{1.5}
\begin{center}
{\scriptsize
\caption{Mean and standard deviation (in braces) of percentage of misclassification and area under the ROC curve for \emph{Golub} and \emph{GTEx} datasets, and those of MSPE, ECP of $50\%$ PI and width of PI for \emph{Eye} dataset for all the competing methods.}
\label{tab2}
\begin{tabular}{|p{1.3 cm}|p{.7 cm}  p{.7 cm}  p{0.7 cm} p{0.7 cm} p{0.7 cm} p{0.7 cm} p{0.7 cm} p{0.7 cm} p{0.7 cm} p{0.7 cm} p{0.7 cm}|}  \hline
Methods$\rightarrow$ \quad Dataset $\downarrow$ &  ISIS-SACD &  ISIS-SACD &  SSR-LASSO &  SSR-EN & SSR-Ridge & \quad PCR &\qquad SPCR & \qquad RPCR & \quad BCR &  RIS-RP & RIS-PCR  \\ \hline
\multicolumn{12}{|c|}{\textbf{\textit{ Misclassification rate and Area under ROC curve for Datasets with Categorical response}}}  \\ \hline
\multicolumn{12}{|c|}{\emph{Mean and SD of Misclassification Rate (in \%)}}  \\ \hline
Golub  & 11.82 (6.90) &  11.50 (7.06) & 45.45 (0.00) & 45.45 (0.00) & 45.45 (0.00) & 7.09 (5.68) & 41.36 (13.31)  & 9.73 (7.28)& 19.36 (9.79)  &  5.54 (4.36) & 5.77 (4.52) \\
GTEx & 0.00 (0.00) & 0.00 (0.00) & 34.83 (0.00) & 34.83 (0.00) & 34.83 (0.00) & 0.06 (0.13) & 3.53 (3.31) & 0.22 (0.18) & 13.28 (3.79) &  0.39 (0.20) & 0.49 (0.32)  \\
\hline
\multicolumn{12}{|c|}{\emph{\qquad Mean and SD of Area under Receiver Operating Characteristic curve}}  \\ \hline
Golub  & 0.876 (.073) &  0.879 (.074) & 0.500 (.000) & 0.500 (.000) & 0.500 (.000) & 0.923 (.062) &0.582 (.134) & 0.895 (.078)&  0.816 (.093) &  0.978 (.027) & 0.943 (.044)  \\
GTEx  & 1.00 (.000) & 1.00 (.000) & 0.500 (.000) & 0.500 (.000) & 0.500 (.000) & 0.999 (.001) & 0.964 (.033) & 0.998 (.001) & 0.877 (.041)  & 1.00 (.000) & .996 (.002)  \\ \hline
\multicolumn{12}{|c|}{\textbf{\textit{ MSPE, ECP of $50\%$ PI} and Width of $50\%$ PI}}  \\ \hline
\multicolumn{12}{|c|}{\emph{Mean and SD of Mean Square Prediction Error}}  \\ \hline
Eye  & 11.66 (4.06) & 11.66 (4.06) & 20.92 (19.33) & 20.92 (19.33) & 7.31 (2.91) & 13.84 (3.94) & 8.65 (3.08) & 7.67 (3.30) & 10.01 (4.04)  & 8.54 (3.09) &  8.29 (2.99)  \\\hline
\multicolumn{12}{|c|}{\emph{Mean and SD of Empirical Coverage Probability and Width of the Prediction Interval}}  \\ \hline
Empirical coverage  & 0.502 (.138) &  0.502 (.138) & 0.634 (.130) & 0.709 (.106) & 0.700 (.076) &  0.423 (.325) & 0.508 (.123) & 0.522 (.114) & 0.564 (.117) &   0.598 (.101) & 0.507 (.107)  \\\hline
Width of interval  & 1.208 (.057) &  1.208 (.057) & 1.970 (.190) & 2.033 (.917) & 1.539 (.303) & 1.884 (1.612) & 1.202 (.079) & 1.055 (.049)  & 1.249 (.056)  & 1.341 (.038) &  1.056 (.036) \\\hline
\end{tabular}}
\end{center}
\end{table*}}

\paragraph{\it Results.} 
For the Golub data set, both the lowest misclassification rate and the highest area under ROC curve are achieved by RIS-RP, which is closely followed by RIS-PCR. TARP based methods attain lower sd than other methods as well. PCR also yields reasonable performance with $7\%$ average misclassification rate and area under ROC more than $0.9$. RPCR and ISIS-based methods produce 
average rates of misclassification of at least $10\%$,  with area under the ROC of $\sim 0.9$. BCR possesses high misclassification rate (about $19\%$), although area under ROC is more than $0.8$. Neither the MSPE, nor the area under ROC curve, is satisfactory for SPCR.  Finally, for all the 100 repetitions, SSR based methods invariably select the intercept-only model. Thus, the MSPEs of these methods depend entirely on the proportion of test samples obtained from the two classes. 

For the GTEx dataset, perfect classification is achieved by ISIS based methods. These methods along with RIS-RP have the highest area under the ROC curve. PCR, RPCR, RIS-RP and RIS-PCR also yield satisfactory results, having less than $0.5\%$ average misclassification rate and more than $99\%$ area under the ROC curve.  SPCR is comparable with an average MSPE of less than $4\%$. BCR attains $13.3\%$ average misclassification rate, with the area under the ROC curve almost $0.9$. SSR based methods fail to show any discriminatory power in the GTEx dataset.  

SSR-Ridge, RPCR, RIS-PCR, RIS-RP and SPCR yield excellent performance in terms of MSPE in the eye dataset with an average MSPE of less than 0.9. SSR-Ridge has an average ECP of about 0.7. RIS-PCR shows more stable performance in terms of ECP, followed by SPCR. BCR and ISIS based methods have similar overall performance. In terms of MSPE, BCR outperforms ISIS based methods but is outperformed by ISIS based methods in terms of ECP. PCR is not quite as good in terms of either measure.  SSR-LASSO and SSR-EN again fail to perform adequately for the eye dataset.

\vskip10pt
\noindent {\bf The GEUVADIS cis-eQTL dataset}~ We conclude this section by illustrating the TARP approach on a massive dataset.
The GEUVADIS cis-eQTL dataset (\cite{massive_p}) is publicly available at \url{http://www.ebi.ac.uk/Tools/geuvadis-das/}. This dataset
consists of messenger RNA and microRNA on lymphoblastoid cell line (LCL) samples from 462 individuals provided by the $1000$ Genomes Project along with roughly $38$ million SNPs. E2F2 plays a key role in the control of the cell cycle. Hence, as in \cite{chen2017}, we choose the gene E2F2 (Ensemble ID: ENSG00000000003) as the response. A total of $8.2$ million $(=p_n)$ SNPs are pre-selected as candidate predictors on the basis of having at least $30$ non-zero expressions. The total number of subjects included in the dataset is about $450$ $(=n)$. The genotype of each SNP is coded as $0$, $1$ or $2$ corresponding to the number of copies of the minor allele.

TARP is applied on this dataset. We consider four different training sample sizes, viz., $n_t=200, 250$, $300$ and $350$, and test sample size $100$ in each case. As $p_n$ is huge, we applied three different values of $\delta$, namely, $2,5$ and $8$, to analyze the effect of a conservative screening. The recommended choice of $\delta$ lies within $(5,6)$ when $p_n=8.2\times 10^6$ and $n\in[200,400]$. To perform SVD for RIS-PCR, we use \emph{fast.svd} instead of the usual \emph{svd} to cope with the massive number of regressors. Table \ref{tab3} provides the MSPE, the ECP of $50\%$ PI and width of the PI, obtained by two different variants of TARP. 
\afterpage{
\begin{table*}[!h]
\renewcommand{\arraystretch}{1.2}
{\scriptsize
\caption{MSPE, ECP and width of PI (in order) obtained by RIS-RP and RIS-PCR for three values of $\delta$ and different training sample sizes $(n_t)$.}
\label{tab3}
\begin{center}
\begin{tabular}{|p{1 cm}|p{2.7 cm}  p{2.7 cm}  p{2.7 cm}| }  \hline
 & \multicolumn{3}{|c|}{RIS-RP}\\
  & $\delta=2$   & $\delta=5$ &  $\delta=8$  \\ 
$n_t $ &MSPE~ECP~ Width & MSPE~ECP~ Width & MSPE~ECP~ Width \\ \hline
$200$ 	& 0.800 \;0.39 \;1.059 & 0.872 \; 0.42 \;0.983 & 0.855 \; 0.34\; 0.928 \\
 $250$ & 0.852\; 0.39\; 1.102 & 0.920 \; 0.42\; 1.023 & 0.921 \; 0.35 \;1.013 \\
 $300$	& 0.860\; 0.36\; 1.126 & 0.855 \; 0.44 \;1.075 &  0.866 \; 0.36 \;1.069 \\
$350$	& 0.778\; 0.45 \;1.210 & 0.779 \; 0.48 \;1.221 & 0.829 \; 0.46 \;1.219 \\\hline
\hline
 & \multicolumn{3}{|c|}{RIS-PCR}\\
  & $\delta=2$   & $\delta=5$ &  $\delta=8$  \\ 
$n_t $ &MSPE~ECP~ Width & MSPE~ECP~ Width & MSPE~ECP~ Width  \\\hline
 200& 0.834 \; 0.06 \; 0.177 & 0.838 \; 0.12 \; 0.192 & 0.831 \; 0.10\; 0.252 \\
 250 &  0.858 \; 0.14 \; 0.355 & 0.882 \; 0.12 \; 0.289 & 0.896 \; 0.19 \;0.420 \\
 300 & 0.845 \; 0.14 \; 0.399 & 0.867 \; 0.20 \; 0.511 & 0.865 \; 0.20 \;0.487 \\
 350 &  0.757 \; 0.35 \; 0.893 & 0.786 \; 0.36 \; 0.886 & 0.826 \; 0.41 \;0.984\\ 
\hline
\end{tabular}
\end{center}}
\end{table*}}
 
 {\it Results:} The MSPEs of RIS-RP and RIS-PCR are comparable for all the choices on $n$. However, RIS-RP yields much better empirical coverage probabilities than RIS-PCR, especially when $n\leq 300$. The three choices of $\delta$ yield comparable results in terms of all the measures in general. For RIS-RP, $\delta=5$ results in higher ECP and for RIS-PCR higher ECP is obtained using $\delta=8$. Moreover, the choice $\delta=8$ makes both the procedures much faster compared to other choices of $\delta$. When the training sample is $350$, $\delta=2,5$ and $8$ select about $290800, 12600$ and $7960$ variables, respectively, on an average in the screening stage out of $8.2\times 10^6$ variables. In view of the results in this massive dimensional dataset, it seems reasonable to use a higher value of $\delta$ for filtering out noisy regressors, and computational convenience.

 \section{Appendix}\label{sec:7}
 
 \begin{small}

This section contains proofs of the theorems stated in the paper.  We use a generic notation $c$ for the constants, although all of them may not be equal.

\noindent {\bf Some Useful Results}

\begin{lemma}\label{lm:1}
 Let  $\varepsilon_n$ be a sequence of positive numbers such that $n\varepsilon_n^2\succ 1$. Then under conditions
 \begin{enumerate}[(a)]
  \item $ln N(\varepsilon_n,\mathcal{P}_n)\leq n\varepsilon_n^2$ for all sufficiently large $n$.
  \item $\pi \left( \mathcal{P}_n^c \right)\leq e^{-2n\varepsilon_n^2}$ for all sufficiently large $n$.
  \item $\pi\left\{ f : d_t(f,f_0) \leq \varepsilon_n^2/4\right\} \geq \exp\{-n\varepsilon_n^2/4\}$ for all sufficiently large $n$ and for some $t>0$.
 \end{enumerate}
Then ~~$P_{f_0}\left[\pi\left\{d(f,f_0)> 4\varepsilon_n| (\yn,{\bf X})\right\} > 2 e^{-n\varepsilon_n^2\left(0.5\wedge (t/4)\right) } \right] \leq 2 e^{-n\varepsilon_n^2\left(0.5\wedge (t/4)\right) }.  $
\end{lemma}
{\noindent The proof is given in \cite{Jiang_2007}.}

\begin{lemma}\label{lm:2}
Suppose assumption (A1) holds. Let $\ad$ be such that $ \sum_j x_j^2 |r_{x_j,y}|^{\delta}/p_n \rightarrow \ad$ as $n\rightarrow\infty$  and ${\bf x}$ be a $p_n\times 1$ sample vector of the regressors. Then the following holds:
\begin{enumerate}[a.]
 \item The random matrix $\Rg$ described in (\ref{eq_rp2}) and (\ref{eq_rp18}) satisfies
\[ \|\Rg {\bf x}\|^2 / p_n\xrightarrow{p} c\ad. \]
\item Let $\|{\bf x}_{\bgm}\|^2=\sum_{j=1}^{p_n} x_j^2 I(\gamma_j=1)$ where $\gamma_j$ is the $j^{th}$ element of the vector $\bgm$ described in (\ref{eq_rp2}) and $I(\cdot)$ is the indicator function, then \[  \| {\bf x}_{\bgm}\|^2 /p_n \xrightarrow{p} c\ad, \]
where $c$ is the proportionality constant in (\ref{eq_rp2}).
 \end{enumerate}
\end{lemma} 
{\noindent The proof is given in the Supplement.}

\begin{lemma}\label{lm:3}
Let $\bta \sim N(0,\sigma_{\theta}^2 I)$, then for a given $\Rg$, ${\bf x}$ and $y$ the following holds
$$P(|(\Rg {\bf  x})^{\prime} \bta - {\bf x}^{\prime} \bbta_0 |< \Delta)> \exp \left\{ - \frac{ ({\bf x}^{\prime} \bbta_0)^2 +\Delta^2 }{\sigma_{\theta}^2 \|\Rg {\bf x}\|^2 } \right\} \frac{2^4 \Delta^4}{\sigma_{\theta}^2 \|\Rg {\bf x}\|^2}.$$ 
\end{lemma}
{\noindent The proof is given in \cite{GD_2015}.}

\vskip5pt
\noindent {\bf Proofs of the Theorems.}
Without loss of generality we consider $|x_j|<M$ with $M=1$, $j=1,2,\ldots,p_n$ and $d(\sigma^2)=1$, although the proofs go through for any fixed value of $M$ and $\sigma^2$.

\begin{proof}[{\bf Proof of Theorem \ref{thm:1}}]

Define the sequence of events 

$\mathcal{B}_n= \left\{ \pi\left\{d(f,f_0)> 4\varepsilon_n| (\yn,{\bf X})\right\} > 2 e^{-n\varepsilon_n^2/4 }  \right\}$ and we need to show $P \left( \mathcal{B}_n^c \right) > 1- 2 e^{-n\varepsilon^2/5}$.
 We first consider the sequence of events $\mathcal{A}_n$ in assumption (A2), and show that $P \left( \mathcal{B}_n^c | \mathcal{A}_n \right) > 1- 2 e^{-n\varepsilon^2/4}$. The proof then follows from assumption (A2) for moderately large $n$.

The proof of $P \left( \mathcal{B}_n^c | \mathcal{A}_n \right) > 1- 2 e^{-n\varepsilon^2/4}$ hinges on showing the three conditions of Lemma \ref{lm:1} for the approximating distribution 
 \begin{eqnarray}
f(y)=\exp\{y a(h)+b(h)+c(y)\}\mbox{ with }h=(\Rg {\bf x})^{\prime} \bta, \label{eq_rp8}  
\end{eqnarray}
 and the true distribution $f_0$, where $\Rg$ is as given in (\ref{eq_rp16}).
 
 \vskip5pt
 \noindent{\it Checking condition (a).}
Let $\mathcal{P}_n$ be the set of densities $f(y)$ stated above with parameter $|\theta_j|<c_n$, $j=1,2,\ldots,m_n$, where $\{c_n\}=\{\sigma_{\theta} \sqrt{5n} \varepsilon_n\}$ and the model $\bgm$ is such that $\bgm \in \mathcal{M}_k$, for some $k\in \{0,1, \ldots, p_n\}$, given $\mathcal{A}_n$. 
For any $\bgm$ the corresponding set of regression parameters can be covered by $\mathit{l}_{\infty}$ balls of the form $B=(v_j-\epsilon,v_j+\epsilon)_{j=1}^{m_n}$ of radius $\epsilon>0$ and center $v_j$. It takes $(c_n/\epsilon+1)^{m_n}$ balls to cover the parameter space for each of model $\bgm$ in $\mathcal{P}_n$. There are at most $\min\{ \binom{p_n}{k}, p_n^{k_n}\}$ models for each $\bgm$ under consideration as we are only concerned with models in $\mathcal{A}_n$ (see assumption (A2)), and there are $(p_n+1)$ possible choices of $k$. Hence it requires at most $ N(\epsilon,k)\leq c p_n^{k_n+1}  (c_n/\epsilon+1)^{m_n}$ $\mathit{l}_\infty$ balls to cover the space of regression parameters $\mathcal{P}_n$, for some constant $c$.  
 
 Next we find the number of Hellinger balls required to cover $\mathcal{P}_n$. We first consider the KL distance between $f$ and $f_0$, then use the fact $d(f,f_0) \leq \left( d_0(f,f_0) \right)^{1/2}$. 
 Given any density in $\mathcal{P}_n$, it can be represented by a set of regression parameters $(u_j)_{j=1}^{m_n}$ falling in one of these $N(\epsilon,k)$ balls $B=(v_j-\epsilon,v_j+\epsilon)_{j=1}^{m_n}$ and $p_{\bgm}=k$. More specifically, let $f_u$ and $f_v$ be two densities in $\mathcal{P}_n$ of the form (\ref{eq_rp8}), where $u=(\Rg {\bf x})^{\prime} \bta_1$, $v=(\Rg {\bf x})^{\prime} \bta_2$ with $|\theta_{i,j}|<c_n$, $i=1,2$ and $p_{\bgm}=k$, then
 \begin{eqnarray}
  d_0(f_u,f_v) &=& \int \int  f_v \log \left( \frac{f_v}{f_u} \right) \nu_y (dy) \nu_{\bf x}(d{\bf x}) \notag \\
  &=& \int \int \left\{ y(a(u)-a(v)) + (b(u)-b(v)) \right\} f_v \nu_y (dy) \nu_{\bf x}(d{\bf x}) \notag \\
 &=& \int (u-v)  \left\{ a^{\prime}(u_v) \left( - \frac{b^{\prime} (v)}{a^{\prime} (v)} \right) +b^{\prime} (u_v) \right\}  \nu_{\bf x}(d{\bf x}). \notag
 \end{eqnarray}
The last expression is achieved by integrating with respect to $y$ and using mean value theorem, where $u_v$ is an intermediate point between $u$ and $v$. Next consider

\hspace{1 in}$ |u-v|= |(\Rg {\bf x})^{\prime} \bta_1 -(R_{\bgm} {\bf x})^{\prime} \bta_2| \leq \| \Rg {\bf x} \| \| \bta_1 - \bta_2 \|, $

using the Cauchy-Schwartz inequality.
Now, by Lemma \ref{lm:2} we have $\|\Rg {\bf x}\|^2/p_n \xrightarrow{p} \ad$ as $n\rightarrow\infty$ for some constant $0<\ad<1$. Therefore we can assume that for sufficiently large $p_n$, $\|\Rg {\bf x}\|\leq   \sqrt{ p_n}$. 
Also, $\|\bta_1-\bta_2 \| \leq \sqrt{m_n} \epsilon$.
 Combining these facts we have $|u-v| \leq  \epsilon \sqrt{m_np_n}$. Similarly $\max\{|u|,|v|\} \leq   c_n \sqrt{m_n p_n} $. 
These together imply that

\vskip5pt
 \hspace{.9 in} $d_0(f_u,f_v) \leq \epsilon \displaystyle\sqrt{m_n p_n}   \left\{ \displaystyle\sup_{|h|\leq  c_n \sqrt{m_n p_n} } |a^{\prime}(h)|   \sup_{|h|\leq  c_n \sqrt{ m_n p_n} } \frac{|b^{\prime}(h)|}{|a^{\prime}(h)|} \right\}. $

\noindent Therefore $d(f_u,f_v) \leq \varepsilon_n$ if we choose 

$\epsilon= \varepsilon_n^2 /\left\{  \sqrt{ m_n p_n} \sup_{|h|\leq   c_n \sqrt{ m_n p_n}} |a^{\prime}(h)|   \sup_{|h|\leq c_n \sqrt{ m_n p_n} } \left( |b^{\prime}(h)| / |a^{\prime}(h)| \right) \right\}$.

Therefore, density $f_u$ falls in a Hellinger ball of size $\varepsilon_n$, centered at $f_v$. As shown earlier, there are at most $N(\epsilon, k)$ such balls. Thus, the Hellinger covering number
\begin{eqnarray}
 N(\varepsilon_n,\mathcal{P}_n ) \leq N(\epsilon, k) =c p_n^{k_n+1}   \left( \frac{c_n}{\epsilon}+1 \right)^{m_n} \notag \hspace{2 in}\\
\hspace{.4 in} = c p_n^{k_n+1} \left[   \left( \frac{c_n}{\varepsilon_n^2} \left\{  \sqrt{m_n p_n} \sup_{|h|\leq  c_n \sqrt{m_n p_n}} |a^{\prime}(h)|   \sup_{|h|\leq  c_n \sqrt{m_n p_n} } \frac{ |b^{\prime}(h)| }{ |a^{\prime}(h)| } \right\} +1 \right) \right]^{m_n} \notag\\
\leq c p_n^{k_n+1} \left( \frac{1}{\varepsilon_n^2}  D( c_n \sqrt{m_n p_n})+1 \right)^{m_n}, \hspace{2.2 in} \notag
\end{eqnarray}
where $D(R)=R \sup_{h\leq R} |a^{\prime} h| \sup_{h\leq R} |b^{\prime}(h)/ a^{\prime} (h)|$. The logarithm of the above quantity is no more than
$\log c + (k_n+1) \log p_n  - m_n \log (\varepsilon_n^2)+ m_n \log \left(1+D(c_n \sqrt{m_n p_n}) \right),$
as $0<\varepsilon_n^2<1$. Using the assumptions in Theorem \ref{thm:1} condition (a) follows.

 \noindent{\it Checking condition (b)} For the $\mathcal{P}_n$ defined in condition (a),
$\pi (\mathcal{P}_n^c) \leq  \pi (\cup_{j=1}^{m_n} |\theta_j|>c_n).$

\noindent Observe that $\pi ( |\theta_j|>c_n) \leq 2 \exp\{-c_n^2/(2\sigma_{\theta}^2 ) \} /\sqrt{2\pi c_n^2/\sigma_{\theta}^2}$ by Mills ratio. Now for the choice that $c_n=\sigma_{\theta} \sqrt{5n} \varepsilon_n$ the above quantity is $2 \exp\{-5 n \varepsilon_n^2/2  \} /\sqrt{10 \pi  n \varepsilon_n^2 }$. Therefore 
$$ \pi (\mathcal{P}_n^c) \leq \sum_{j=1}^{m_n} \pi ( |\theta_j|>c_n) \leq 2 m_n\exp\{-5 n \varepsilon_n^2 /2\}/\sqrt{10 \pi  n \varepsilon_n^2 } \leq e^{-2n\varepsilon_n^2}$$
for sufficiently large $n$. Thus condition (b) follows.  

 \noindent{\it Checking condition (c)}
Condition (c) is verified for $t=1$. Observe that
$$ d_{t=1}(f,f_0)=\int \int f_0 \left( \frac{f_0}{f} -1 \right) \nu_y (dy) \nu_{\bf x}(d{\bf x}). $$
Integrating out $y$ we would get $\int E_{y|{\bf x}} \left[ \left\{ (f_0/f) (Y) -1 \right\} \right] \nu_{\bf x}(d{\bf x}). $ Note that under $f$ and $f_0$ we have same function of $y$ as given in (\ref{eq_rp8}) with $h={\bf x}^{\prime} \bbta_0$ for $f_0$. Therefore, the above can be written as $E_{\bf x} \left[ \left\{ (\Rg {\bf x})^{\prime} \bta - {\bf x}^{\prime} \bbta_0 \right\} g\left(u^{*} \right)  \right]$ using mean value theorem where $g$ is a continuous derivative function, and $u^{*}$ is an intermediate point between $(\Rg {\bf x})^{\prime} \bta$ and ${\bf x}^{\prime} \bbta_0$. Therefore, if $\left| (\Rg {\bf x})^{\prime} \bta - {\bf x}^{\prime} \bbta_0 \right|< \Delta_n$, then $|u^{*}|< |{\bf x}^{\prime} \bbta_0 |+\Delta_n$. This in turn implies that for sufficiently small $\Delta_n$, $|g(u^{*})|$ will be bounded, say by $M$.
Consider a positive constant $\Delta_n$. From Lemma \ref{lm:3} we have
\begin{eqnarray}
 P(| (\Rg {\bf x})^{\prime} \bta - {\bf x}^{\prime} \bbta_0|< \Delta_n ) &=& \sum_{\bgm} P(| (\Rg {\bf x})^{\prime} \bta - {\bf x}^{\prime} \bbta_0|< \Delta_n | \bgm ) \pi(\bgm) \notag \\
 &\geq & E_{\bgm} \left[ \exp \left\{ - \frac{ ({\bf x}^{\prime} \bbta_0)^2 +\Delta_n^2 }{\sigma_{\theta}^2 \|\Rg {\bf x}\|^2 } \right\} \frac{2^4 \Delta^4}{\sigma_{\theta}^2 \|\Rg {\bf x}\|^2} \right] \notag \\
 &=& \frac{2^4 \Delta_n^4}{ ({\bf x}^{\prime} \bbta_0)^2 +\Delta_n^2 } E_{\bgm} \left\{ \frac{Z_{\bgm}}{p_n} \exp \left(-\frac{Z_{\bgm}}{p_n} \right)   \right\}, \label{eq_rp11}
\end{eqnarray}
where $Z_{\bgm}=\left\{ ({\bf x}^{\prime} \bbta_0)^2 +\Delta_n^2 \right\}/ \left\{\sigma_{\theta}^2 \|\Rg {\bf x} \|^2 / p_n \right\}$. By part (a) of Lemma \ref{lm:2}, and continuous mapping theorem $Z_{\bgm} \xrightarrow{p} \left\{ ({\bf x}^{\prime} \bbta_0)^2 +\Delta_n^2 \right\}/ \left(\sigma_{\theta}^2 c \ad \right) > \Delta_n^2 / \left(\sigma_{\theta}^2 c\ad \right)$.

For some non-negative random variable $Z$ and non-random positive numbers $p$, $a$ and $b$, consider the following fact
 \vspace{-.1 in}
\begin{eqnarray}
 E\left(\frac{Z}{p} \exp \left\{ -\frac{Z}{p} \right\} \right) &\geq& a P\left(\frac{Z}{p} \exp \left\{ -\frac{Z}{p} \right\} > a \right) \notag \\
 &\geq& a P\left(\frac{Z}{p}> \frac{a}{b}, \exp \left\{ -\frac{Z}{p}\right\} >ab \right) \notag \\
 &=& a P \left(Z>\frac{ap}{b}, Z< - p\log (ab)\right)  \notag \\
  &=& a P\left(\frac{ap}{b} < Z< - p\log (ab)\right). \label{eq_rp10}
\end{eqnarray}
Replacing $Z$ by $Z_{\bgm}$, $p$ by $p_n$ and taking $a=\Delta_n^2 \exp \{-n\varepsilon_n^2/3 \}/ ( \sigma_{\theta}^2 c\ad)$, and $b=p_n$ $ \exp \{-n\varepsilon_n^2/3 \}$. Thus 
$-p_n \log (ab) = -p_n $ $ \log \left[ \Delta_n^2 p_n \exp \{-2n\varepsilon_n^2/3 \} / (\sigma_{\theta}^2 c \ad ) \right]> p_n n \varepsilon_n^2/2$ and  $ap_n/b= \Delta_n^2/ \left(\sigma_{\theta}^2 c \ad \right)$ 
for sufficiently large $n$. Therefore the expression in (\ref{eq_rp10}) is greater than 
$$\frac{\Delta_n^2}{\sigma_{\theta}^2 c \ad} e^{-n\varepsilon_n^2/3 }  P\left(\frac{\Delta_n^2}{\sigma_{\theta}^2 c\ad}  \leq Z_{\bgm} \leq \frac{1}{2} p_n n \varepsilon_n^2 \right)  $$
Note that $({\bf x}^{\prime} \bbta_0)^2 <\sum_{j=1}^{p_n} |\beta_{0,j}|<K$, and the probability involved in the above expression can be shown to be bigger than some positive constant $p$ for sufficiently large $n$. Using these facts along with equation (\ref{eq_rp11}), we have $ P(| (\Rg {\bf x})^{\prime} \bta - {\bf x}^{\prime} \bbta_0|< \Delta_n ) > \exp\{ -n\varepsilon^2/4 \} $. Choosing $\Delta_n <\varepsilon^2/(4M)$ condition (c) follows. 
\end{proof}

 
 \begin{proof}[{\bf Proof of Theorem \ref{thm:2}}]
 The outline of the proof of Theorem \ref{thm:2} closely follows the arguments given in the proof of Theorem \ref{thm:1}. Therefore we only present those parts of the proof which are different. As in Theorem \ref{thm:1}, we show that $P \left( \mathcal{B}_n^c | \mathcal{A}_n \right) > 1- 2 e^{-n\varepsilon^2/4}$ by checking the three conditions of Lemma \ref{lm:1}. 
 
 The proof of Condition (a) is the same as for Theorem \ref{thm:1}, except for the places involving the projection matrix $\Rg$. Observe that given a dataset $\mathcal{D}^n$ and other tuning parameters we fix a particular projection matrix $\Rg$. The only property of $\Rg$ needed to prove condition (a) is $\|\Rg {\bf x} \|^2   \leq p_n $ for sufficiently large $n$. To show this we use that fact that $\Rg$ is a matrix with orthonormal row vectors, and therefore $\|\Rg {\bf x} \|^2 \leq \| {\bf x}_{\bgm} \|^2  \leq p_n $.
 
 The proof of Condition (b) depends only on the prior assigned on $\bta$, and therefore remains the same under the conditions of Theorem \ref{thm:2}.
 
 The proof of Condition (c) differs from that of Theorem \ref{thm:1} in showing $ P(| (\Rg {\bf x})^{\prime} \bta - {\bf x}^{\prime} \bbta_0|< \Delta_n ) > \exp\{ -n\varepsilon^2/4 \} $ for some constant $\Delta_n$. To see this consider a positive constant $\Delta_n$. As before, from Lemma \ref{lm:3} we have
\vspace{-.1 in}
 \begin{eqnarray}
 P(| (\Rg {\bf x})^{\prime} \bta - {\bf x}^{\prime} \bbta_0|< \Delta_n ) 
 & \geq  E_{\bgm} \left[ \exp \left\{ - \displaystyle\frac{ ({\bf x}^{\prime} \bbta_0)^2 +\Delta_n^2 }{\sigma_{\theta}^2 \|\Rg {\bf x}\|^2 } \right\} \displaystyle\frac{2^4 \Delta^4}{\sigma_{\theta}^2 \|\Rg {\bf x}\|^2} \right] \notag \\
& \geq  E_{\bgm} \left[ \exp \left\{ - \displaystyle\frac{ ({\bf x}^{\prime} \bbta_0)^2 +\Delta_n^2 }{\sigma_{\theta}^2 \alpha_n \| {\bf x}_{\bgm}\|^2 } \right\} \displaystyle\frac{2^4 \Delta^4}{\sigma_{\theta}^2 \| {\bf x}_{\bgm}\|^2} \right] \notag \\
 & = \displaystyle\frac{2^4 \Delta_n^4}{ ({\bf x}^{\prime} \bbta_0)^2 +\Delta_n^2 } E_{\bgm} \left\{ \displaystyle\frac{Z_{\bgm}}{p_n} \exp \left(-\frac{Z_{\bgm}}{\alpha_n p_n} \right)   \right\}, \label{eq_rp19}
\end{eqnarray}
where $Z_{\bgm}=\left\{ ({\bf x}^{\prime} \bbta_0)^2 +\Delta_n^2 \right\}/ \left\{ \sigma_{\theta}^2 \| {\bf x}_{\bgm} \|^2 / p_n \right\}$, and $\alpha_n$ is as in (A3). From part (b) of Lemma \ref{lm:2}, and continuous mapping theorem we have $Z_{\bgm} \xrightarrow{p} \left\{ ({\bf x}^{\prime} \bbta_0)^2 +\Delta_n^2 \right\}/\left( \sigma_{\theta}^2 c\ad \right)$ $> \Delta_n^2 / \left(\sigma_{\theta}^2 c \ad \right)$.

For some positive random variable $Z$ and non-random positive numbers $p$, $a$ and $b$, consider the following 
\vspace{-.2 in}
\begin{eqnarray}
 E\left(\frac{Z}{p} \exp \left\{ -\frac{Z}{\alpha p} \right\} \right) &\geq& a P\left(\frac{Z}{p} \exp \left\{ -\frac{Z}{\alpha p} \right\}\geq a \right) \\
 &\geq& a P\left(\frac{Z}{p} \geq \frac{a}{b}, \exp \left\{ -\frac{Z}{\alpha p}\right\}  \geq ab \right) \notag  \\
 &=& a P\left(\frac{ap}{b} < Z< - \alpha p\log (ab)\right). \label{eq_rp20}
\end{eqnarray}
Replacing $Z$ by $Z_{\bgm}$, $p$ by $p_n$, $\alpha$ by $\alpha_n$ and taking $a=\Delta_n^2 \exp \{-n\varepsilon_n^2/3 \}/ (\sigma_{\theta}^2 c \ad)$, and $b=p_n \exp \{-n\varepsilon_n^2/3 \}$. Thus, 
\begin{eqnarray*}
-\alpha_n p_n \log (ab) = -\alpha_n p_n \log \left[ \Delta_n^2 p_n \exp \{-2n\varepsilon_n^2/3 \} / ( \sigma_{\theta}^2 c \ad ) \right]
\sim 2 p_n \log \left( \Delta_n^2 p_n / ( \sigma_{\theta}^2 c \ad ) \right) /3 > 2 p_n /3
\end{eqnarray*}
 for sufficiently large $n$ and $ap_n/b= \Delta_n^2/ \left( \sigma_{\theta}^2 c \ad \right). $
Therefore the expression in (\ref{eq_rp20}) is greater than 

\vskip5pt
\hspace{1.5 in}$\displaystyle\frac{\Delta_n^2}{\sigma_{\theta}^2 c \ad} e^{-n\varepsilon_n^2/3 }  P\left(\frac{\Delta_n^2}{\sigma_{\theta}^2 c \ad}  \leq Z_{\bgm} \leq \frac{2}{3} p_n \right).  $

\vskip5pt
Note that $({\bf x}^{\prime} \bbta_0)^2 <\sum_{j=1}^{p_n} |\beta_{0,j}|<K$, and the probability involved in the above expression can be shown to be bigger than some positive constant $p$ for sufficiently large $n$. Using these facts along with equation (\ref{eq_rp19}), we have $ P(| (\Rg {\bf x})^{\prime} \bta - {\bf x}^{\prime} \bbta_0|< \Delta_n ) > \exp\{ -n\varepsilon^2/4 \} $. Choosing $\Delta_n <\varepsilon^2/(4M)$ condition (c) follows. 
 \end{proof}
\bibliographystyle{apalike}
 \bibliography{rp.bib}
\end{small}
 \end{document}


\maketitle

\section{Additional Simulation Results}

\subsection{Results for Larger $p_n$}
Here we present the performance of the different methods with respect to mean square prediction error (MSPE) and empirical coverage probability (ECP) of $50\%$ prediction intervals (PI), along with the width of PIs for a larger choice of $p_n$ in each scheme.

\afterpage{
\begin{figure}
\begin{center}
     \includegraphics[scale=.3]{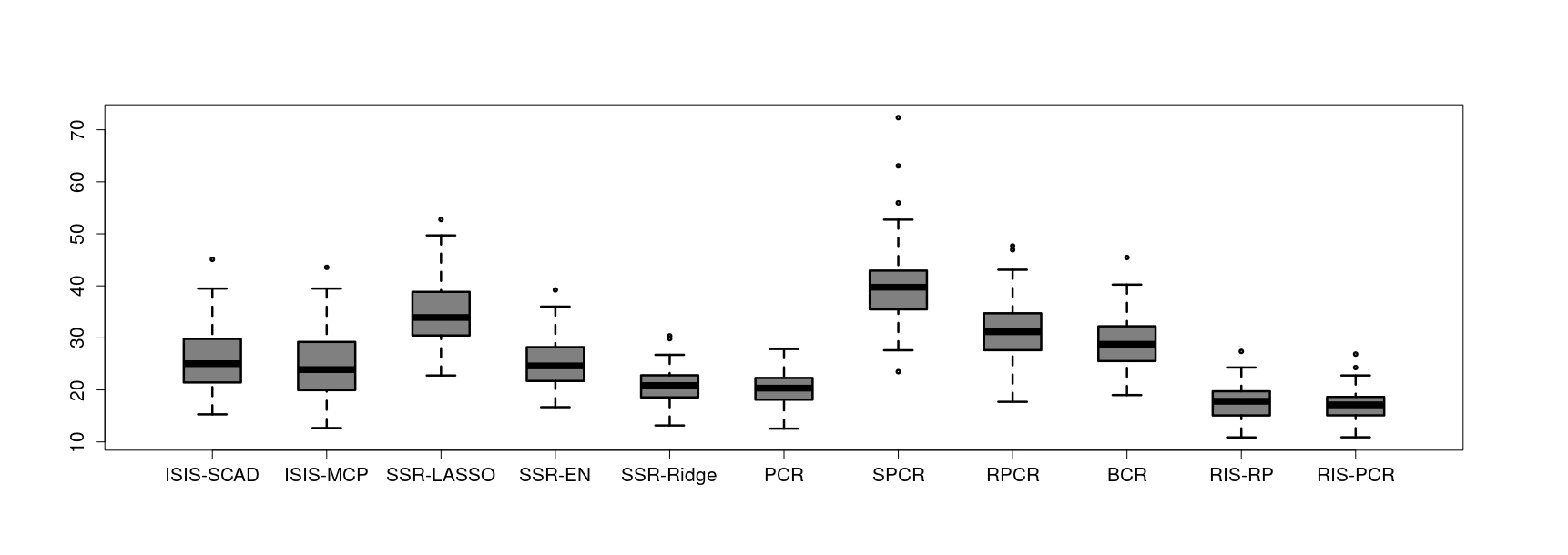}
\caption{Box-plot of MSPE of eleven competing methods for $p_n=3000$ in \emph{Scheme I}. }
\label{fig2}
\end{center}
\end{figure}
 }

 \paragraph{Mean Square Prediction Error (MSPE).} In \emph{Scheme I} (see Figure \ref{fig2}) for $p_n=3000$, TARP based methods yield the best results, with RIS-PCR having marginally better performance than RIS-RP. These methods are followed by SSR-Ridge and PCR. ISIS based methods and SSR-EN have overall comparable performance with SSR-EN showing more stable results, and ISIS-MCP having minimum average MSPE among these 3 methods. Performance of these six methods are followed by BCR, and that in turn is followed by RPCR. SSR-LASSO and SPCR have very poor performance in this scenario. 
 
 As in the main paper simulation case, TARP based methods have the best overall performance in \emph{Scheme II} for $p_n=2*10^4$, immediately followed by PCR (see Figure \ref{fig4}). These three methods are followed by SSR-EN and SSR-Ridge. Unlike for $p_n=10^4$, here SSR-EN outperforms SSR-Ridge. RPCR is substantially worse but is still better than SSR-LASSO, BCR and SPCR (in that order). ISIS based methods had very poor performance.
 
In \emph{Scheme III} (see Figure \ref{fig6}), the performance of all the methods are almost similar, except SSR based methods. Among the SSR based methods, SSR-EN and SSR-LASSO have very similar results. SSR-EN has lower average MSPE, however SSR-LASSO has a more stable performance. Unlike for $p_n=5000$, ISIS based methods often yield large values of MSPE for $p_n=10^4$. 
 
The relative performance of SPCR and ISIS based methods are unsatisfactory in \emph{Scheme IV} (see Figure \ref{fig9}).
Among the other methods (see Figure \ref{fig10}), TARP based methods and PCR have the best performance, followed by BCR and SSR-EN. The relative performance of RPCR and SSR-LASSO improve significantly for $p_n=2*10^4$. These two methods outperform SSR-Ridge, and among them, RPCR yields much better results than SSR-LASSO.

 \afterpage{
\begin{figure}
\begin{center}
\includegraphics[scale=.3]{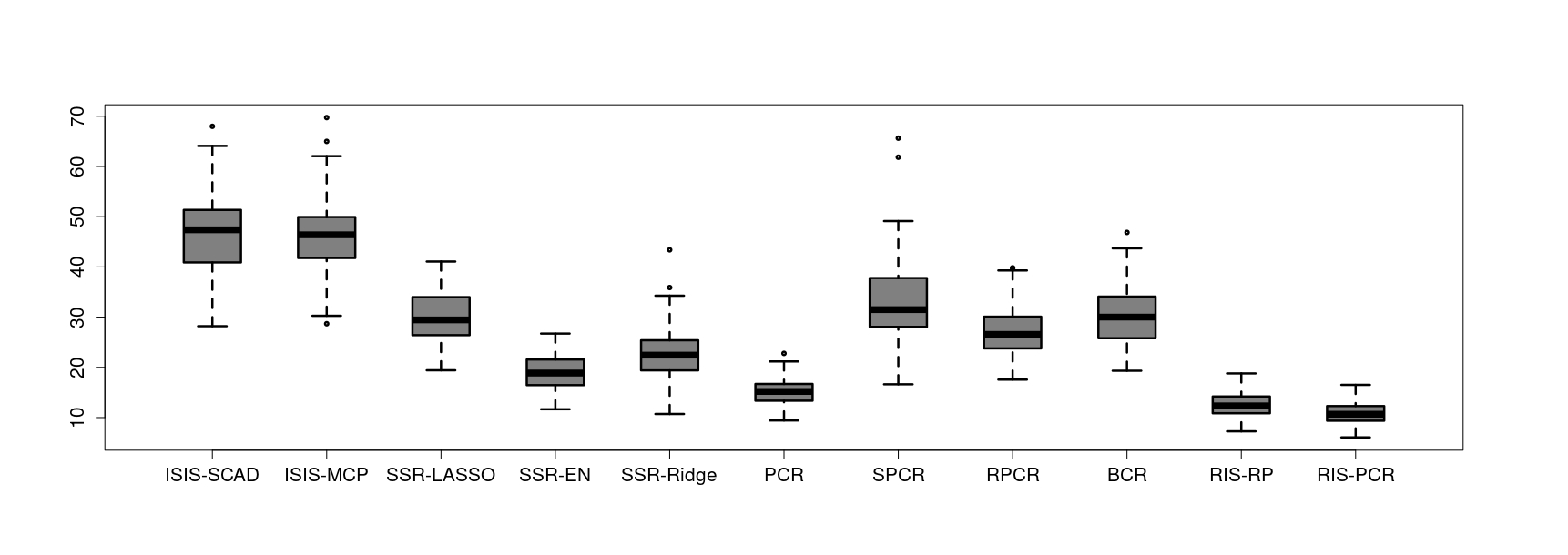}
\caption{Box-plot of MSPE of eleven competing methods for $p_n=20000$ in \emph{Scheme II}. }
\label{fig4}
\end{center}
\end{figure}
 }
 \afterpage{
\begin{figure}
\begin{center}
\includegraphics[scale=.3]{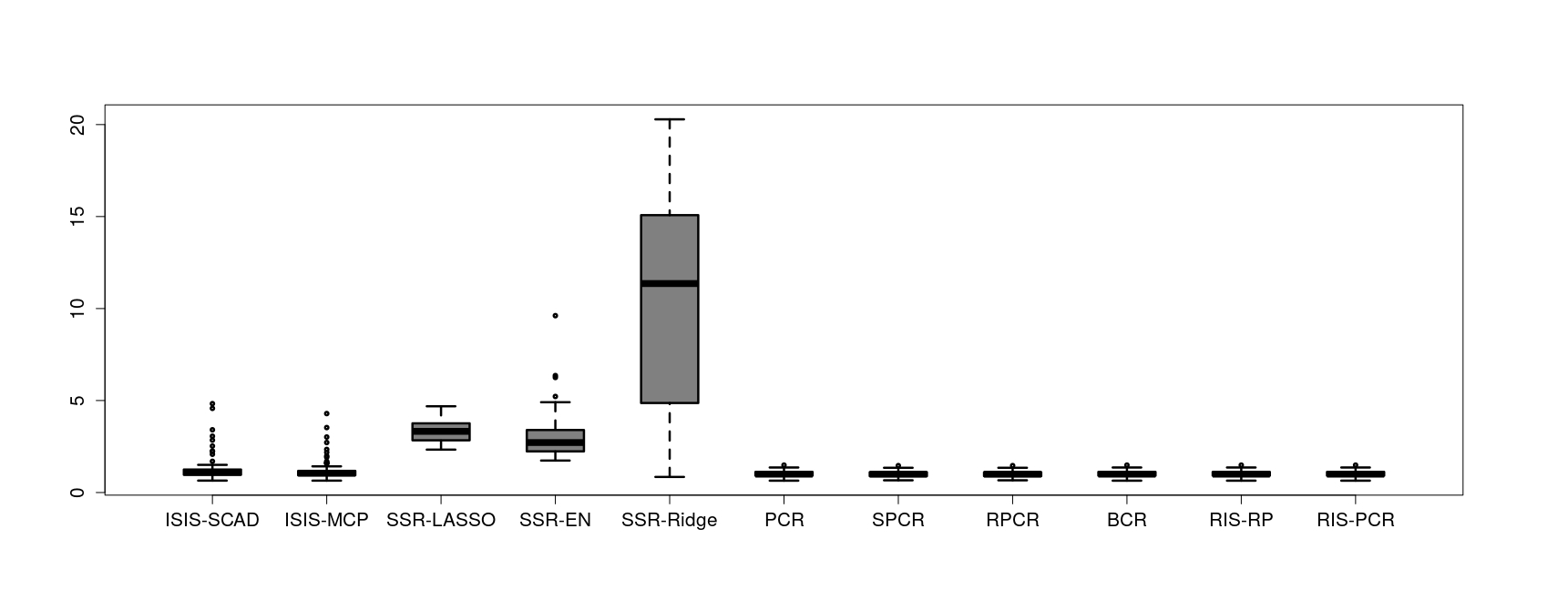}
\caption{Box-plot of MSPE of eleven competing methods for $p_n=10000$ in \emph{Scheme III}. }
\label{fig6}
\end{center}
\end{figure}
 }

\afterpage{
\begin{figure}
\centering     
\subfigure[MSPE of all the methods.]{\label{fig9}  \includegraphics[scale=.3]{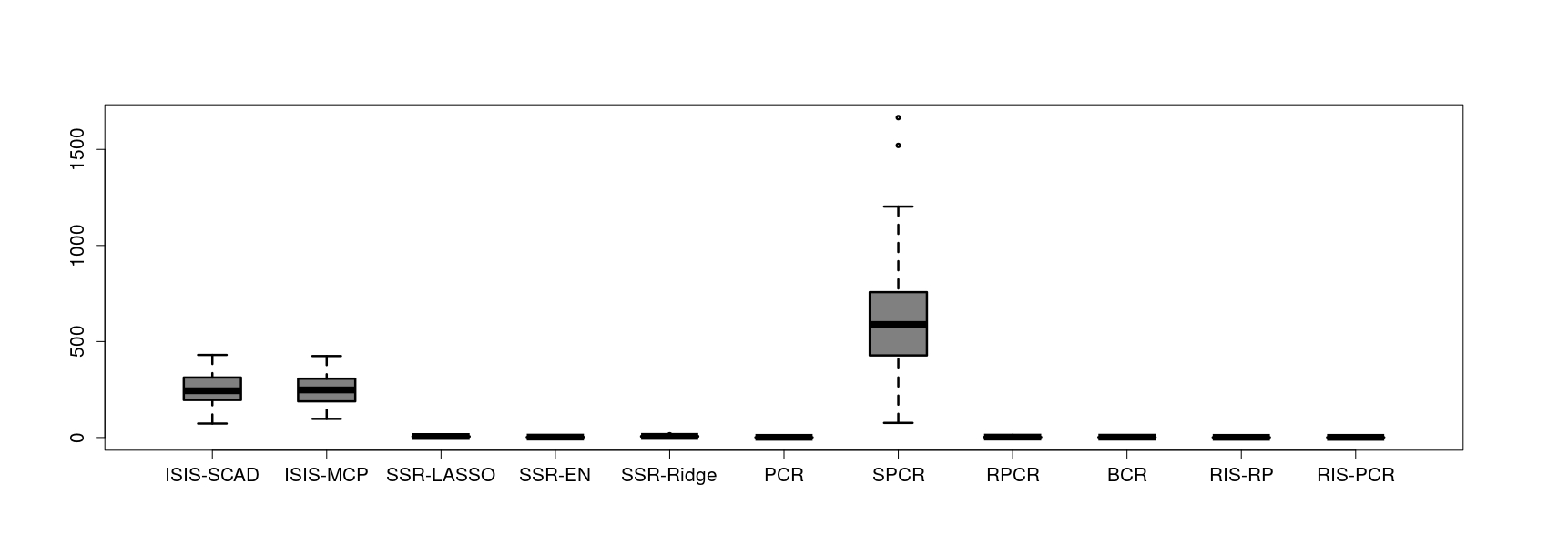}}
\subfigure[MSPE of selected methods.]{\label{fig10}  \includegraphics[scale=.3]{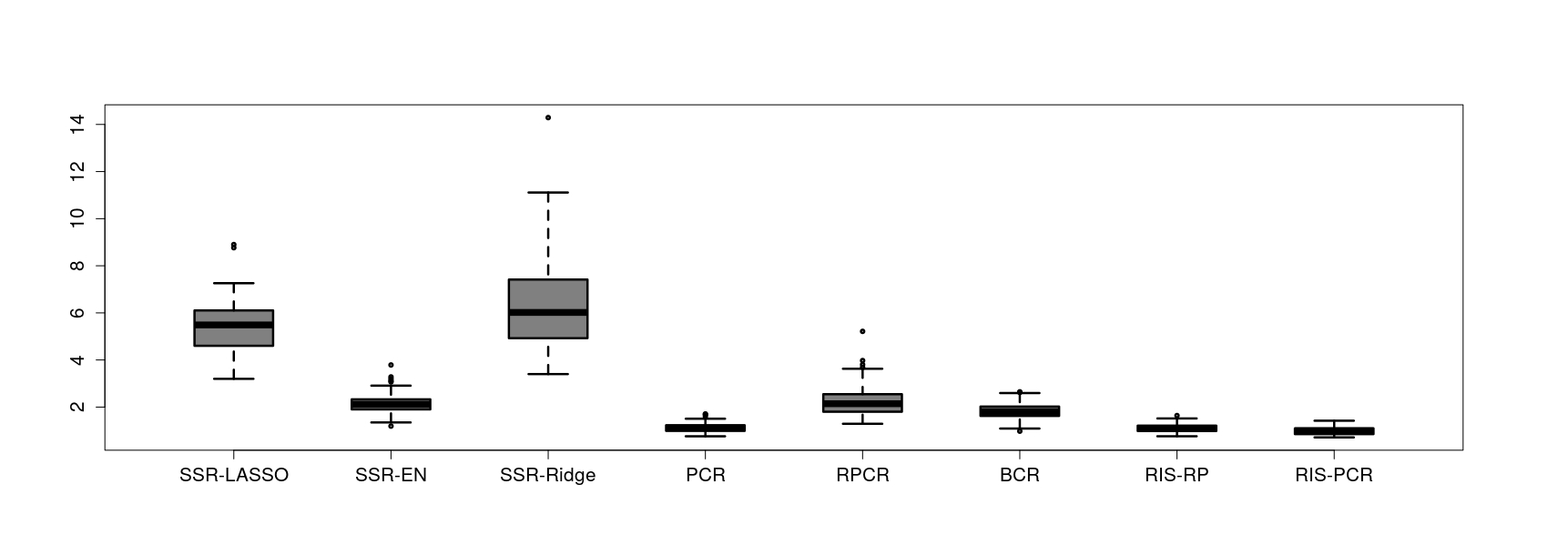}}
\caption{Box-plot of MSPE of the competing methods for $p_n=20000$ in \emph{Scheme IV}. } 
\end{figure}
 }

 \afterpage{
\begin{table*}[!h]
\renewcommand\thetable{4}
\renewcommand{\arraystretch}{1.5}
\begin{center}
{\scriptsize
\caption{Mean and standard deviation(sd) of empirical coverage probabilities of 50\% prediction intervals, and mean and sd of width of 50\% prediction intervals.}
\label{tab4}
\begin{tabular}{|p{1.6 cm}|p{0.65 cm} p{0.65 cm} p{0.65 cm} p{0.65 cm} p{0.65 cm} p{0.65 cm} p{0.65 cm} p{0.65 cm} p{0.65 cm} p{0.65 cm} p{0.75 cm}|}  \hline
Methods$\rightarrow$ \quad Scheme,~ $p_n$  &  ISIS-SACD &  ISIS-SACD &  SSR-LASSO &  SSR-EN &  SSR-Ridge & \qquad PCR & \qquad SPCR & \qquad RPCR & \qquad BCR & RIS-RP & RIS-PCR \\ \hline 
\multicolumn{12}{|c|}{\emph{Average and standard deviation (in braces) of empirical coverage probability}} \\ \hline
I~~ \;~$3*10^3$ & 0.307 (.052) & 0.306 (.050)  & 0.434 (.059)  & 0.536 (.043)  & 0.631 (.048)  & 0.437 (.288)  & 0.498 (.059) & 0.471 (.055)  & 0.283 (.053)  & 0.464 (.052)  & 0.382 (.054)  \\ \hline

II~  \; ~$2*10^4$ & 0.334 (.051) & 0.333 (.048) & 0.422 (.055) & 0.577 (.045) & 0.637 (.049) & 0.365 (.260) & 0.505 (.059) & 0.454 (.052) &  0.280 (.058) &  0.488 (.056) & 0.312 (.055)   \\ \hline

III \; ~$10^4$ & 0.495 (.053) & 0.494 (.053) &  0.502 (.058) & 0.653 (.033) &  0.671 (.144) & 0.503 (.058) & 0.498 (.052)  & 0.498 (.052) &  0.550 (.054) & 0.487 (.055) &   0.499 (.054)  \\ \hline

IV~  \; ~$2*10^4$ & 0.432 (.063) &  0.433 (0.060)  & 0.485 (.099)  &  0.658 (.034)  & 0.664 (.108) & 0.489 (.308) & 0.492 (.052) & 0.483 (.058) & 0.494 (.051)  & 0.675 (.060) &  0.408 (.052)  \\ 
\hline
\multicolumn{12}{|c|}{\emph{Average and standard deviation (in braces) of width of the 50\% prediction interval}} \\ \hline
I~~   \;~$3*10^3$ &  4.096 (.522) &  3.942 (.531) & 6.148 (.410) & 7.388 (.543) & 8.098 (.484) & 8.082 (7.354) & 8.441 (.608) & 7.023 (.492) & 5.720 (.411) & 5.192 (.152) & 4.090 (.335)  \\ \hline

II~  \; ~$2*10^4$ & 5.964 (.554) & 5.954 (.554) & 5.276 (.525) & 7.039 (.479) & 8.604 (.674) & 6.589 (7.380)  & 7.649 (.929) & 6.270 (.565) &  5.687 (.465) & 4.594 (.204)  & 2.635 (.209) \\ \hline

III  \; ~$10^4$ & 1.483 (.174) &  1.452 (.144) & 1.36 (.072) & 1.942 (.122) & 5.439 (2.018) & 1.362 (.072) & 1.373 (.070) & 1.373 (.070) & 1.551 (.073) & 1.343 (.079) & 1.383 (.077)  \\ \hline

IV~ \; ~$2*10^4$ & 18.00 (2.876) & 17.98 (2.789) & 3.428 (2.522) & 1.621 (.167) & 4.087 (.634) & 2.895 (3.808)  & 32.430 (7.647) & 1.925 (.234) & 1.794 (.144) & 2.066 (.142) &  1.070 (.071)\\ 
\hline
\end{tabular}
}
\end{center}
\end{table*}
 }

 \paragraph{Empirical Coverage Probability (ECP) and Prediction Interval (PI) Widths.}
 The relative performance of all the methods is quite similar to that observed in the main paper (see Table 1 of the paper and Table \ref{tab4}) for all the simulation schemes. The detailed description is given 
below. 
 
The average ECP of SPCR is closest to $0.5$ in \emph{Scheme I}, followed by RPCR, RIS-RP and SSR-EN, with RIS-RP showing the smallest average width among these methods. The sd of ECPs and widths of the PIs is highest for PCR, indicating a lack of stability in it's performance. RIS-PCR yields some under coverage, however the average width of PI for RIS-PCR is lowest among all the competing methods. SSR-Ridge has an average ECP of about $0.63$ and a high width of the PI. ISIS based methods and BCR have lower-coverage probabilities on average.

For \emph{Scheme II}, SSR-RP and SPCR have the best overall results, with SPCR having higher PI widths than SSR-RP. The ECPs of SSR-LASSO, SSR-EN and RPCR are closer to $0.5$ than the remaining methods with SSR-LASSO having a lower average width of the PIs. The average ECP of PCR is close to $0.5$, however there are huge fluctuations in the widths of PIs (sd$=7.4$). Among the other methods, SSR-Ridge shows higher coverage with a comparatively large width of PI. The remaining four methods show under coverage for both the choices of $p_n$. Among these four methods, RIS-PCR attains the lowest average width of PI.

All the methods, except SSR-EN and SSR-Ridge have comparable performance in terms of both the criteria in \emph{Scheme III}. SSR-EN and SSR-Ridge tend to have higher ECP and larger width of PI in this setup.  

SSR-LASSO, RPCR and BCR yield the best performance in \emph{Scheme IV}, and BCR achieves the lowest average width of PI among them. Although the average ECP of SPCR is close to $0.5$, it possesses much higher width than all other methods. 
SSR-Ridge, SSR-EN and RIS-RP show higher ECP in general. In terms of PI widths of these three methods, performance of RIS-RP is the best, followed by SSR-EN.
In \emph{Scheme IV} as well, PCR has a higher sd of the width of PIs indicating instability in its performance. ISIS-SCAD, ISIS-MCP and RIS-PCR show some under coverage, and averages of ECPs for all of them are around $0.4$. RIS-RCR yields the lowest overall width of PI among all the methods. Average width of PI for ISIS based methods and SPCR are much higher than all the other methods in \emph{Scheme IV}. 
 

 
\subsection{Effect of sample size on different competing methods} 
We consider three different choices of $n$, viz., $150, 200, 250$. The results for $n=200$ are presented in the paper and previous sections of Supplementary file. The results corresponding two other choices of $n$ are presented in this section. The relative performance of the competing methods for different values of $n$ remains similar for the two choices of $p_n$. Therefore, we present results corresponding to higher values of $p_n$ only. 

 \paragraph{Comparison with respect to MSPE.}
 
\afterpage{
\begin{figure}
\centering     
\subfigure[$n=150$]{\label{fig11}  \includegraphics[scale=.275]{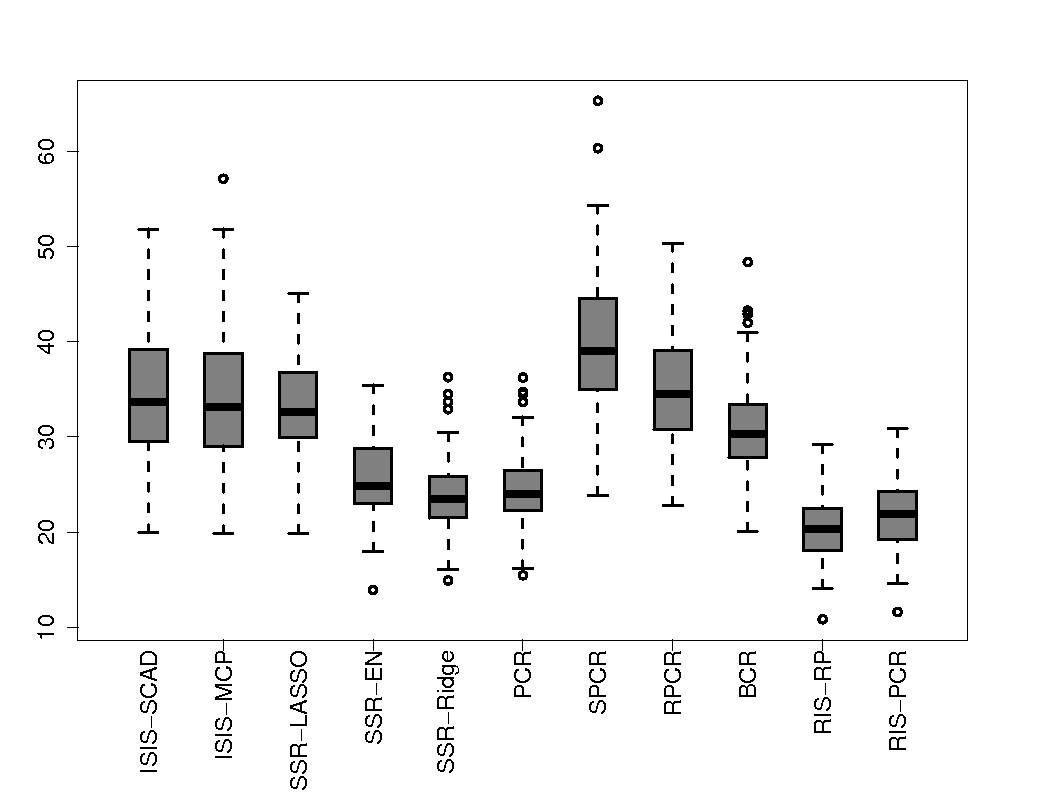}}
\subfigure[$n=250$]{\label{fig12}  \includegraphics[scale=.275]{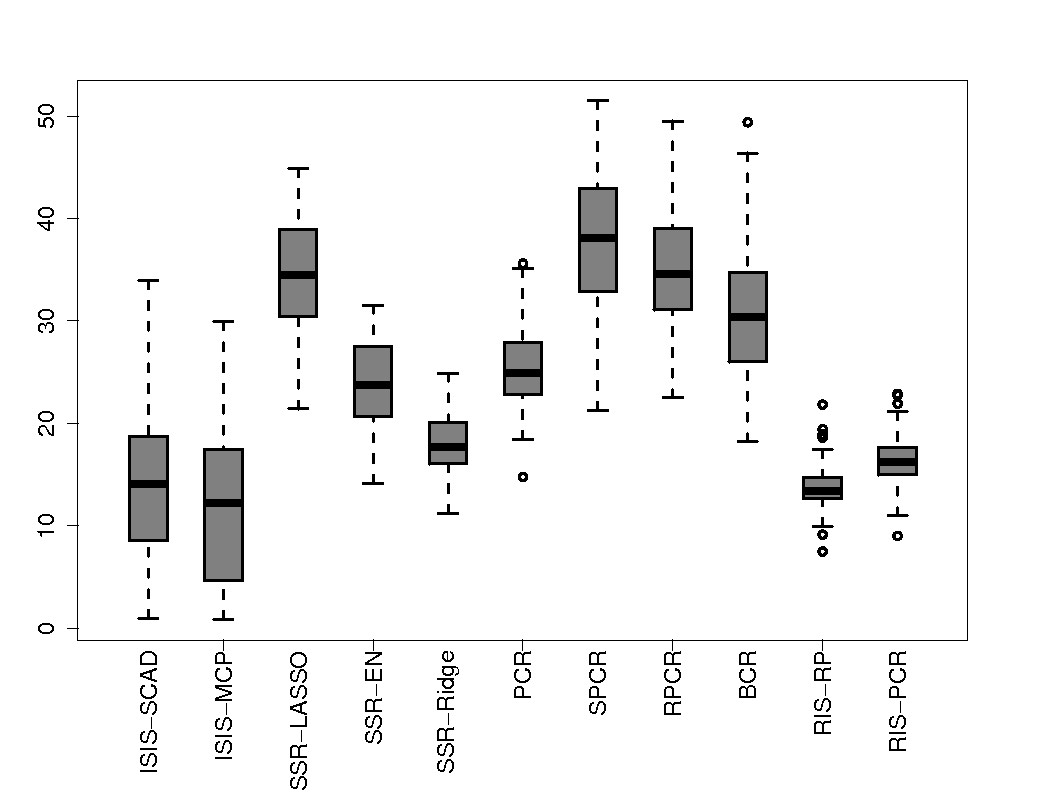}}
\caption{Box-plot of MSPE of the competing methods for $p_n=3\times 10^3$ in \emph{Scheme I}. } 
\end{figure}
 }

 In \emph{Scheme I} (see Figures \ref{fig11}-\ref{fig12}), SSR-Ridge, SSR-EN, PCR and TARP based methods work well for lower sample size. Among these methods, RIS-RP has the best overall performance, immediately followed by RIS-PCR. Among the other methods, BCR and SSR-LASSO yield better performance than others. The overall scenario changes reasonably for higher sample size, when $n=250$. ISIS-MCP based methods outperform other methods with respect to average MSPE, closely followed by RIS-RP, ISIS-SCAD and RIS-PCR which also have satisfactory performance. However, the ISIS based methods have much larger dispersion than TARP based methods, indicating lack of stability in their results. The other methods do not show much improvement with increment of sample size. Among these methods, SSR-Ridge yields relatively better results followed by SSR-EN and PCR.

 \afterpage{
\begin{figure}
\centering     
\subfigure[$n=150$]{\label{fig13}  \includegraphics[scale=.275]{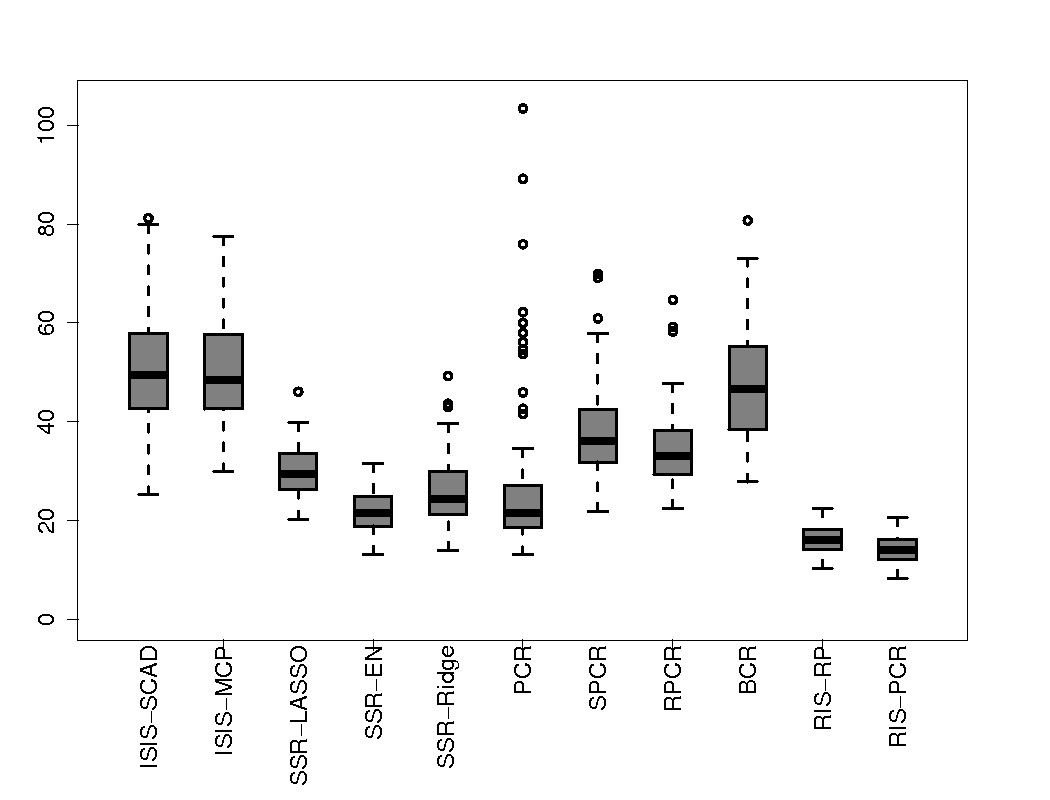}}
\subfigure[$n=250$]{\label{fig14}  \includegraphics[scale=.275]{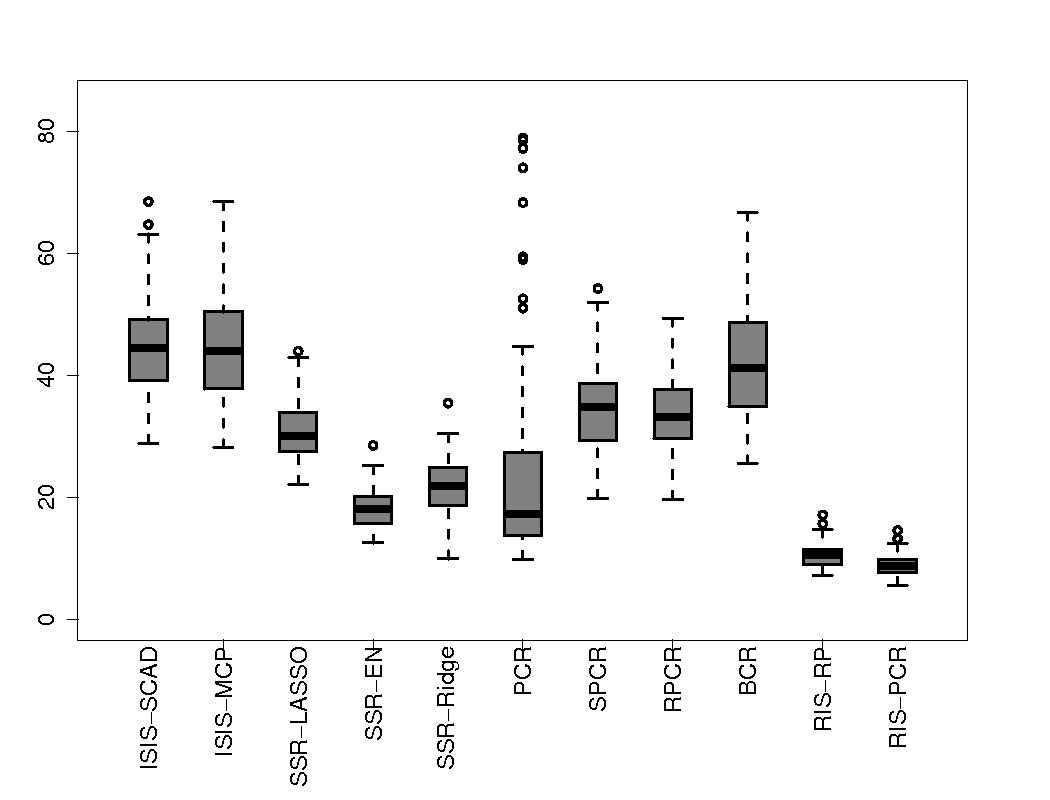}}
\caption{Box-plot of MSPE of the competing methods for $p_n=2\times 10^4$ in \emph{Scheme II}. } 
\end{figure}
 }

 The relative performance of all these methods does not vary much for different sample sizes in \emph{Scheme II} (see Figures \ref{fig13}-\ref{fig14}). For both the choices of $n$, TARP outperforms all other methods. SSR based methods show reasonable results, although the differences of their performance with TARP increases with increment of sample size. The other methods do not perform well for either choices of $n$. Notably, PCR yields a low median MSPE in general, however, it results in very high values of MSPE occasionally. For $n=150$, the MSPE is as large as $1038$, while for $n=250$ the MSPE increases up to $3089$.
 
\afterpage{
\begin{figure}
\centering     
\subfigure[$n=150$]{\label{fig15}  \includegraphics[scale=.275]{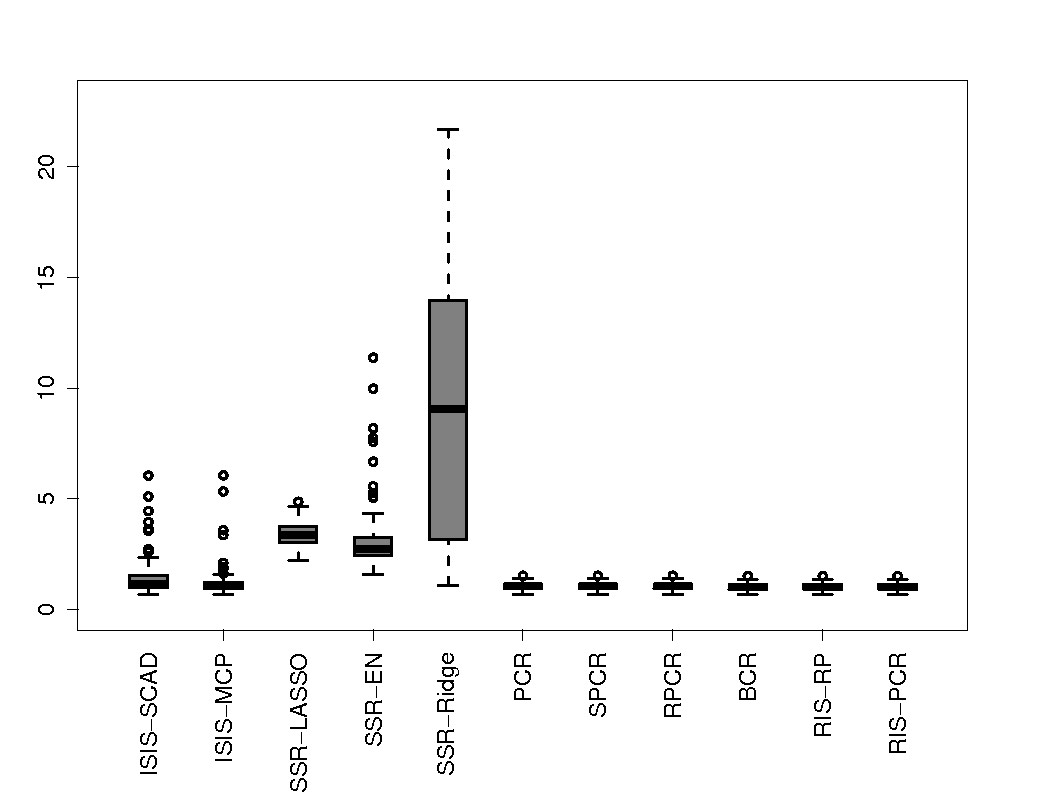}}
\subfigure[$n=250$]{\label{fig16}  \includegraphics[scale=.275]{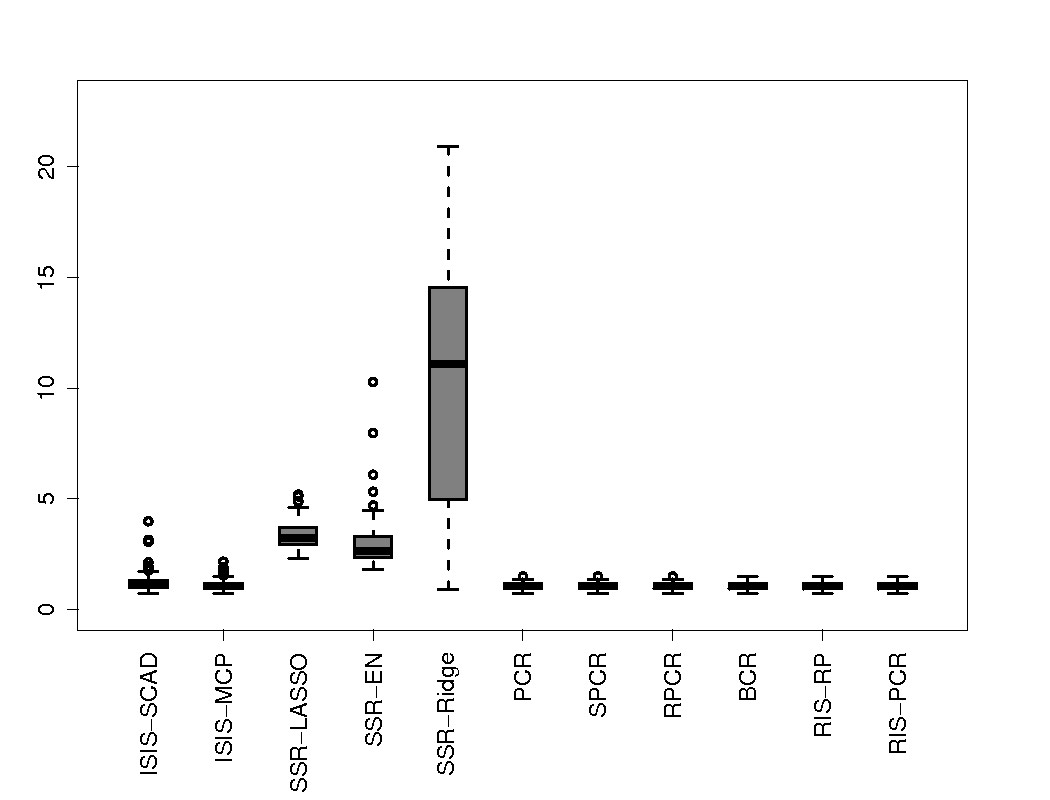}}
\caption{Box-plot of MSPE of the competing methods for $p_n=10^4$ in \emph{Scheme III}. } 
\end{figure}
}

Projection based methods work extremely well in \emph{Scheme III} for both choices of $n$ (see Figures \ref{fig15}-\ref{fig16}). SSR based methods fail to perform well in this scheme, with SSR-Ridge showing largest overall MSPE. ISIS based methods show low MSPE in general, however they occasionally yield higher values of MSPEs.

\afterpage{
\begin{figure}
\centering     
\subfigure[$n=150$]{\label{fig17}  \includegraphics[scale=.275]{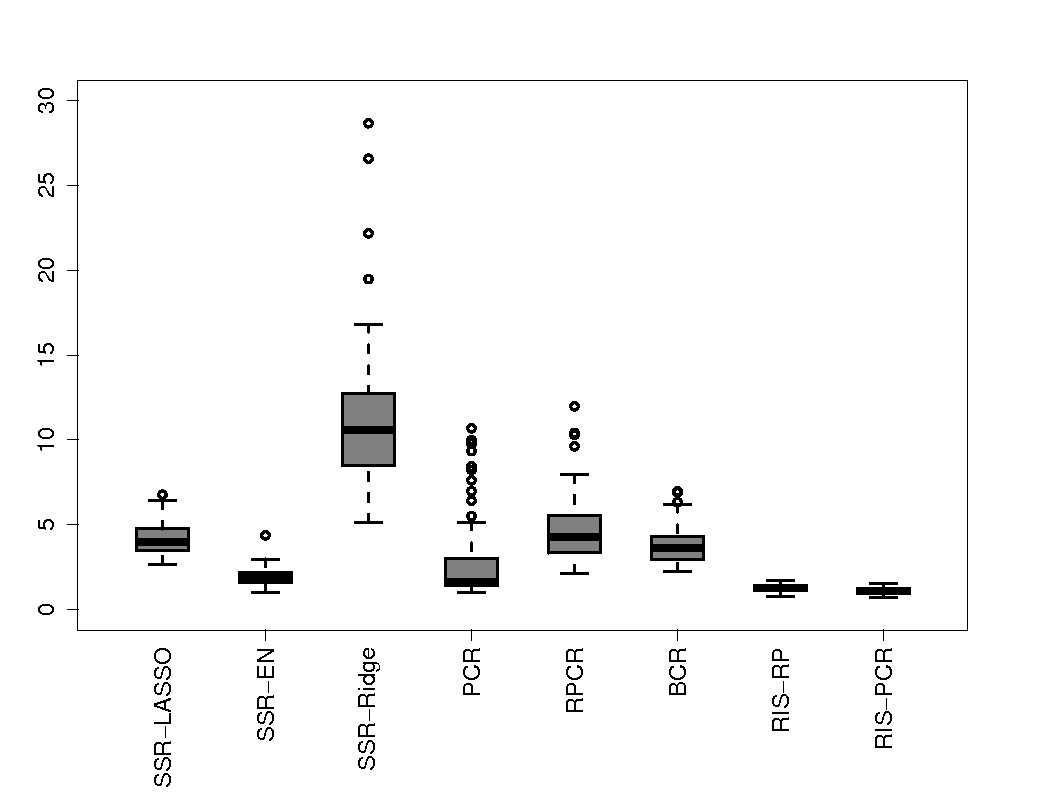}}
\subfigure[$n=250$]{\label{fig18}  \includegraphics[scale=.275]{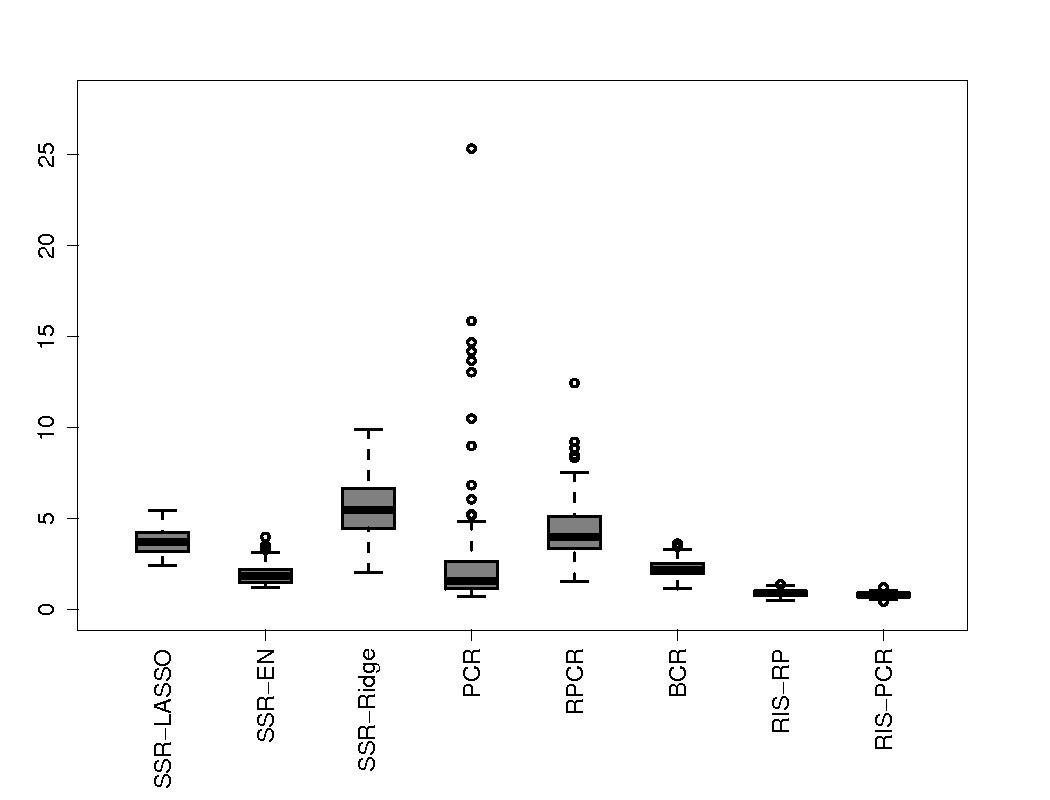}}
\caption{Box-plot of MSPE of the competing methods for $p_n=2\times 10^4$ in \emph{Scheme IV}. } 
\end{figure}
}

As the ISIS based method and SPCR yield extremely large values of MSPE in \emph{Scheme IV}, we avoid presenting further results on those methods (see Figures \ref{fig17}-\ref{fig18}). TARP yields the best performance for both choices of $n$ in this scheme. SSR-EN also gives comparable performance. Among the other methods, SSR-LASSO and BCR yield reasonable results. The other three methods tend to have higher values of MSPE in this scheme. As in \emph{Scheme II}, PCR tends to have extremely large MSPEs. For $n=250$, the highest MSPE obtained by PCR is around $125$, whereas the average MSPE of all other methods is below $7$.

 \afterpage{
\begin{table*}[!h]
\renewcommand\thetable{5}
\renewcommand{\arraystretch}{1.5}
\begin{center}
{\scriptsize
\caption{Mean and standard deviation (in small font) of empirical coverage probabilities of 50\% prediction intervals.}
\label{tab5}
\begin{tabular}{|p{2.2 cm}|p{1.25 cm} p{1.25 cm} p{1.25 cm} p{1.25 cm}|p{1.25 cm} p{1.25 cm} p{1.25 cm} p{1.25 cm} |}  \hline
$n$  $\rightarrow$   &  \multicolumn{4}{c|}{$n=150$} &  \multicolumn{4}{c|}{$n=250$} \\ 
 Schemes $\rightarrow$    &  I  &  II  &  III  & VI  &  I  &  II  &  III  & VI  \\ 
  Methods$\downarrow$ $p_n\rightarrow$   &  $3\times 10^3$ &  $2\times 10^4$ &  $10^4$ & $2\times 10^4$ &  $3\times 10^3$ & $2\times 10^4$ & $ 10^4$ & $2\times 10^4$ \\ \hline 

ISIS-SCAD & $0.290_{\small .060}$ & $0.340_{\small .068}$ & $0.488_{\small .064}$ & $0.429_{\small .057}$ & $0.344_{\small .060}$ & $0.349_{\small .057}$ & $0.507_{\small .061}$ & $0.447_{\small .050}$\\

ISIS-MCP  & $0.294_{\small .052}$ & $0.336_{\small .070}$ & $0.489_{\small .063}$ & $0.427_{\small .058}$ & $0.345_{\small .065}$ & $0.354_{\small .059}$ & $0.503_{\small .061}$ & $0.443_{\small .053}$ \\ 

SSR-LASSO & $0.406_{\small .059}$ & $0.386_{\small .053}$ & $0.494_{\small .064}$ & $0.660_{\small .031}$  & $0.441_{\small .052}$ & $0.436_{\small .050}$ & $0.502_{\small .061}$ & $0.656_{\small .049}$ \\ 

ISIS-EN   & $0.658_{\small .035}$ & $0.660_{\small .035}$ & $0.655_{\small .031}$ & $0.657_{\small .030}$ & $0.656_{\small .031}$ &  $0.661_{\small .035}$ & $0.657_{\small .033}$ & $0.656_{\small .030}$ \\ 

ISIS-Ridge& $0.660_{\small .031}$ & $0.652_{\small .034}$ & $0.652_{\small .032}$ & $0.656_{\small .030}$ & $0.661_{\small .031}$ & $0.660_{\small .034}$ & $0.659_{\small .032}$ & $0.657_{\small .033}$ \\ 

PCR       & $0.407_{\small .200}$ & $0.450_{\small .285}$ & $0.494_{\small .064}$ & $0.484_{\small .287}$ & $0.421_{\small .193}$ & $0.416_{\small .275}$ & $0.502_{\small .061}$ & $0.504_{\small .316}$   \\ 

SPCR      & $0.507_{\small .064}$ & $0.499_{\small .055}$ & $0.494_{\small .064}$ & $0.490_{\small .056}$ & $0.491_{\small .054}$ & $0.497_{\small .050}$ & $0.502_{\small .061}$ & $0.501_{\small .054}$ \\ 

RPCR      & $0.487_{\small .052}$ & $0.475_{\small .056}$ & $0.494_{\small .065}$ & $0.490_{\small .053}$ & $0.489_{\small .057}$ & $0.484_{\small .051}$ & $0.502_{\small .061}$ & $0.494_{\small .058}$ \\ 

BCR       & $0.519_{\small .056}$ & $0.353_{\small .062}$ & $0.626_{\small .061}$ & $0.525_{\small .058}$ & $0.555_{\small .057}$ & $0.353_{\small .054}$ & $0.557_{\small .059}$ & $0.510_{\small .052}$ \\ 

RIS-RP    & $0.508_{\small .063}$ & $0.460_{\small .059}$ & $0.492_{\small .065}$ & $0.723_{\small .050}$ & $0.598_{\small .058}$ & $0.504_{\small .057}$ & $0.499_{\small .059}$ & $0.704_{\small .047}$ \\

RIS-PCR   & $0.434_{\small .062}$ & $0.270_{\small .049}$ & $0.508_{\small .066}$ & $0.395_{\small .058}$ & $0.549_{\small .054}$ & $0.338_{\small .051}$ & $0.510_{\small .060}$ & $0.420_{\small .052}$ \\ 

\hline
\end{tabular}
}
\end{center}
\end{table*}
 }

 \paragraph{Empirical Coverage Probability (ECP) and Prediction Interval (PI) Widths.} Tables \ref{tab5} and \ref{tab6} provide the means and standard deviations (sd) of empirical coverage probabilities (ECP) of $50\%$ prediction intervals (PI) and the widths of the PIs, respectively, for two different choices of $n$.
 
 As $n$ increases, ISIS based methods tend to improve rapidly with respect to both the measures. The other methods show little difference with an increment of $n$. ISIS based methods show under coverage for all the simulation schemes except \emph{Scheme III}. However, the PI widths are low in general for these methods, except for \emph{Scheme IV}, where they result in the highest average width. SSR-EN and SSR-Ridge tend to have higher coverage probabilities, and higher average PI width. For most of the schemes, these two methods possess the highest average width. PCR tends to show under coverage for \emph{Schemes I} and \emph{II}, and yield reasonable ECPs for the other two schemes. However, the widths of PIs for PCR attain the highest variance in all the schemes, except \emph{Scheme III}, indicating lack of stability in the results. SPCR and RPCR yield the best performance in terms of ECP. The width of PI is higher on an average for SPCR than RPCR. For \emph{Scheme IV}, SPCR yields extremely large widths of PI. Except for \emph{Scheme II}, BCR and SSR-LASSO perform reasonably well in terms of ECP. However, these two methods result in some over-coverage in \emph{Scheme III} and \emph{Scheme IV}. With respect to PI width, SSR-LASSO outperforms BCR. TARP based methods yield the lowest PI width in general. RIS-RP yields an average ECP close to $0.5$ in most of the cases, and tend to show high coverage probabilities occasionally. However, RIS-PCR tends to show low coverage probabilities in general, especially for $n=150$.

 \afterpage{
\begin{table*}[!h]
\renewcommand\thetable{6}
\renewcommand{\arraystretch}{1.5}
\begin{center}
{\scriptsize
\caption{Mean and standard deviation (in small font) of width of 50\% prediction intervals.}
\label{tab6}
\begin{tabular}{|p{2.2 cm}|p{1.25 cm} p{1.25 cm} p{1.25 cm} p{1.25 cm}|p{1.25 cm} p{1.25 cm} p{1.25 cm} p{1.25 cm} |}  \hline
$n$  $\rightarrow$   &  \multicolumn{4}{c|}{$n=150$} &  \multicolumn{4}{c|}{$n=250$} \\ 
 Schemes $\rightarrow$    &  I  &  II  &  III  & VI  &  I  &  II  &  III  & VI  \\ 
  Methods$\downarrow$ $p_n\rightarrow$   &  $3\times 10^3$ &  $2\times 10^4$ &  $10^4$ & $2\times 10^4$ &  $3\times 10^3$ & $2\times 10^4$ & $ 10^4$ & $2\times 10^4$ \\ \hline 

ISIS-SCAD & $4.39_{\small 0.50}$ & $6.13_{\small 0.85}$ & $7.87_{\small 44.45}$ & $13.76_{\small 2.60}$ & $2.99_{\small 0.87}$ & $5.87_{\small 0.54}$ & $1.49_{\small 0.25}$ & $9.67_{\small 4.11}$\\

ISIS-MCP  & $4.41_{\small 0.54}$ & $6.12_{\small 0.86}$ & $1.46_{\small 0.30}$ & $14.20_{\small 2.24}$ & $2.78_{\small 0.95}$ & 
$5.90_{\small 0.55}$ & $1.40_{\small 0.14}$ & $9.71_{\small 4.01}$ \\ 

SSR-LASSO & $5.80_{\small 0.41}$ & $4.95_{\small 0.38}$ & $1.36_{\small 0.08}$ & $3.92_{\small 0.40}$  & $6.45_{\small 0.42}$ & $5.49_{\small 0.50}$ & $1.36_{\small 0.06}$ & $3.76_{\small 0.47}$ \\ 

ISIS-EN   & $9.66_{\small 0.77}$ & $8.87_{\small 0.79}$ & $3.35_{\small 0.69}$ & $2.65_{\small 0.34}$ & $9.35_{\small 0.78}$ &  $8.11_{\small 0.74}$ & $3.27_{\small 0.53}$ & $2.67_{\small 0.38}$ \\ 

ISIS-Ridge& $9.31_{\small 0.75}$ & $9.64_{\small 1.19}$ & $5.39_{\small 2.05}$ & $6.32_{\small 0.99}$ & $8.12_{\small 0.64}$ & 
$8.83_{\small 0.95}$ & $5.73_{\small 1.97}$ & $4.49_{\small 0.63}$ \\ 

PCR       & $5.64_{\small 2.89}$ & $7.97_{\small 8.33}$ & $1.36_{\small 0.08}$ & $2.60_{\small 2.82}$ & $5.91_{\small 2.97}$ & 
$6.77_{\small 6.41}$ & $1.36_{\small 0.06}$ & $2.82_{\small 3.53}$   \\ 

SPCR      & $8.44_{\small 0.61}$ & $8.22_{\small 0.80}$ & $1.36_{\small 0.08}$ & $31.08_{\small 9.79}$ & $8.31_{\small 0.53}$ & 
$7.93_{\small 0.70}$ & $1.36_{\small 0.06}$ & $30.21_{\small 10.41}$ \\ 

RPCR      & $7.76_{\small 0.50}$ & $7.73_{\small 0.64}$ & $1.36_{\small 0.08}$ & $2.82_{\small 0.49}$ & $7.72_{\small 0.45}$ & 
$7.54_{\small 0.54}$ & $1.36_{\small 0.06}$ & $2.75_{\small 0.46}$ \\ 

BCR       & $7.44_{\small 0.48}$ & $6.28_{\small 0.52}$ & $1.81_{\small 0.07}$ & $2.74_{\small 0.23}$ & $8.18_{\small 0.42}$ & 
$5.97_{\small 0.38}$ & $1.55_{\small 0.05}$ & $2.05_{\small 0.14}$ \\ 

RIS-RP    & $6.12_{\small 0.26}$ & $4.91_{\small 0.22}$ & $1.35_{\small 0.08}$ & $2.46_{\small 0.16}$ & $6.26_{\small 0.21}$ & 
$4.35_{\small 0.17}$ & $1.36_{\small 0.06}$ & $1.97_{\small 0.13}$ \\

RIS-PCR   & $5.81_{\small 0.33}$ & $2.53_{\small 0.20}$ & $1.40_{\small .08}$ & $1.09_{\small 0.08}$ & $6.10_{\small 0.37}$ & $2.63_{\small 0.20}$ & $1.39_{\small 0.06}$ & $0.98_{\small 0.06}$ \\ 

\hline
\end{tabular}
}
\end{center}
\end{table*}
 }

\section{Predictive Calibration in Binary Response Datasets}
 Apart from measuring the misclassification rates and the area under ROC curve, we validate TARP in terms of it's ability to quantify uncertainly in real datasets with binary responses. To this end, we partition the interval $[0,1]$ into ten equal sub-intervals, viz., $[0,0.1), ~[0.1,0.2)$ and so on, and classify the test data points $({\bf x},y)_{i,new}$ to the $k^{th}$ class if predictive probability of $y_{i,new}$ falls in that class. Next, we consider the squared difference of the empirical proportion of $y_{i,new}=1$ among the data points classified in a given interval with the middle point of the interval, and consider the mean of these squared differences (\emph{MSD}) of all the intervals. If a method is well calibrated, then the \emph{MSD} would be small. The following table shows means and standard deviations of \emph{MSD}s of eleven competing methods for Golub and GTEx datasets.

 \afterpage{
\begin{table*}[!h]
\renewcommand\thetable{7}
\renewcommand{\arraystretch}{1.5}
\begin{center}
{\scriptsize
\caption{Mean and standard deviation (in bracket) of mean square differences (\emph{MSD}) of empirical and predictive probabilities.}
\label{tab7}
\begin{tabular}{|p{1.3 cm}|p{.7 cm}  p{.7 cm}  p{0.7 cm} p{0.7 cm} p{0.7 cm} p{0.7 cm} p{0.7 cm} p{0.7 cm} p{0.7 cm} p{0.7 cm} p{0.7 cm}|}  \hline
Methods$\rightarrow$ \quad Dataset $\downarrow$ &  ISIS-SACD &  ISIS-SACD &  SSR-LASSO &  SSR-EN & SSR-Ridge & \quad PCR &\qquad SPCR & \qquad RPCR & \quad BCR &  RIS-RP & RIS-PCR  \\ \hline
Golub  & 2.605 (.036) &  2.606 (.035) & 2.804 (.000)  & 2.804 (.000) & 2.742 (.028) & 2.589 (.012) & 3.429 (.073)  & 2.587 (.017) & 2.886 (.032)  &  2.611 (.045) & 2.555 (.044) \\
GTEx & 2.784 (.000) & 2.784 (.000) & 3.325 (.000) & 3.325 (.000) & 3.325 (.000) & 2.784 (.000) & 3.216 (.102) & 2.784 (.000) & 2.873 (.007) &  2.782 (.001) & 2.783 (.001)  \\

\hline
\end{tabular}
}
\end{center}
\end{table*}
 }

Table \ref{tab7} indicates that the competing methods do not show much differences in MSD of empirical and predictive probabilities. ISIS based methods, PCR, RPCR and TARP perform relatively well compared to the other methods. Among these methods, RIS-PCR and RIS-RP have lowest MSD for Golub and GTEx, respectively. Relative performance of SSR based methods are deficient in terms of MSD as well. SPCR results in high MSD for both the datasets too.

\section{Mathematical Details}
\begin{proof}[Proof of Lemma 2 a.]
 Consider the conditional expectation and variance of $\|\Rg {\bf x}\|^2$ given $(\bgm,{\bf x})$ as follows:
 \begin{eqnarray*}
   E\left( \|\Rg {\bf x}\|^2|\bgm \right)&=&\|{\bf x}_{\gamma}\|^2 \\
   var\left( \|\Rg {\bf x}\|^2|\bgm \right)&=& 4 \|{\bf x}_{\gamma}\|^2 \left\{1+\left(\frac{1}{\psi}-1 \right) \frac{\sum_{j=1}^{\pg} x_{\gamma,j}^4}{2 \|{ \bf x_{\gamma}}\|^2} \right\},
  \end{eqnarray*}
where ${\bf x}_{\bgm}$ includes the regressors $j$ for which $\gamma_j=1$. The detailed proof is given in Result \ref{res1} below.

Next consider the conditional expectation of $\|\Rg {\bf x}\|^2$ given ${\bf x}$ is given by
\begin{eqnarray}
 E_{\bgm}   E\left( \|\Rg {\bf x}\|^2|\bgm \right)=E_{\bgm} \left(\sum_j x_j^2 I(\gamma_j=1) \right)=c \sum_j x_j^2 |r_{{\bf x}_j,\yn}|^{\delta}, \label{eq_rp4}
\end{eqnarray}
where $c>0$ is the proportionality constant. 
 Also the conditional variance of $\|\Rg {\bf x}\|^2$ given ${\bf x}$ is given by
\begin{eqnarray}
  var_{\bgm} \left\{ E\left( \|\Rg {\bf x}\|^2|\bgm \right) \right\} + E_{\bgm} \left\{  var\left( \|\Rg {\bf x}\|^2|\bgm \right)\right\}. \label{eq_rp3}
\end{eqnarray}
Considering both the terms in (\ref{eq_rp3}) separately as follows:
\begin{eqnarray}
  var_{\bgm} \left\{ E\left( \|\Rg {\bf x}\|^2|\bgm \right) \right\} &=& var_{\bgm} \left( \sum_j x_j^2 I(\gamma_j=1)\right) \notag \\
  &=& c \sum_j x_j^4 |r_{{\bf x}_j,\yn}|^{\delta} \left( 1 - c|r_{{\bf x}_j,\yn}|^{\delta} \right) \leq p_n, \label{eq_rp6}
\end{eqnarray}
as given ${\bf x}$, $\gamma_j$s are independent, and each $|x_j|\leq 1$, and $q_j= c|r_j|^{\delta} <1$. Again
\begin{eqnarray}
 E_{\bgm} \left\{  var\left( \|\Rg {\bf x}\|^2|\bgm \right)\right\} &=& E_{\bgm} \left[4 \|{\bf x}_{\gamma}\|^2 \left\{1+\left(\frac{1}{\psi}-1 \right) \frac{\sum_{j=1}^{\pg} x_{\gamma,j}^4}{2 \|{\bf x_{\gamma}}\|^2} \right\} \right] \notag \\
 &\leq &  c E_{\bgm} \left[  \|{\bf x}_{\gamma}\|^2    \right] \notag \\
 &\leq & c \sum_j x_j^2 |r_{{\bf x}_j,\yn}|^{\delta}, \label{eq_rp7}
\end{eqnarray}
for some constant $c$, as $\sum_{j=1}^{\pg} x_{\gamma,j}^4 < \|{\bf x_{\gamma}}\|^2$.

Therefore, from (\ref{eq_rp4}), (\ref{eq_rp6}) and (\ref{eq_rp7}) it can be shown that the expectation of $ \|\Rg {\bf x}\|^2 /p_n$ converges to the limit $c\ad$, and variance of the same converges to $0$.

\vskip15pt
{\it Proof of Lemma 2 b.}
The proof follows from observing that $$E_{\bgm} \left( \|{\bf x}_{\bgm}\|^2 \right)=c \sum_j x_j^2 |r_{{\bf x}_j,\yn}|^{\delta}$$ and 
$$ var_{\bgm} \left\{ E\left( \|{\bf x}_{\bgm}\|^2|\bgm \right) \right\}c \sum_j x_j^4 |r_{{\bf x}_j,\yn}|^{\delta} \left( 1 - c|r_{{\bf x}_j,\yn}|^{\delta} \right) \leq p_n.$$
Therefore it can be shown that the expectation of $ \| {\bf x}_{\bgm}\|^2 /p_n$ converges to the limit $c\ad$, and variance of the same converges to $0$.
\end{proof}


\begin{result}
\label{res1} 
 Consider a random matrix $\Rg$ which depends on another random vector $\bgm$ distributed as in (2). Then the conditional distribution of $\Rg$ satisfies the following:
 \begin{enumerate}[a.]
  \item $ E\left( \|\Rg {\bf x}\|^2|\bgm \right)=\|{\bf x}_{\gamma}\|^2$, and
\item    $var\left( \|\Rg {\bf x}\|^2|\bgm \right) = 4   \|{\bf x}_{\gamma}\|^2 \left\{1+\left( \psi^{-1}-1\right) \sum_{j=1}^{\pg} x_{\gamma,j}^4 /\left(2\|{\bf x_{\gamma}}\|^2 \right) \right\} .$
 \end{enumerate}
\begin{proof}[Proof of part a.] Observe that 
\begin{eqnarray}
\|\Rg {\bf x}\|^2 &=& \left\| \left( \sum_j r_{1,j} \gamma_j x_j, \sum_j r_{2,j} \gamma_j x_j, \ldots , \sum_j r_{m_n,j} \gamma_j x_j \right)^{\prime} \right\|^2 \notag \\
&=& \left( \sum_j r_{1,j} \gamma_j x_j \right)^2 + \left( \sum_j r_{2,j} \gamma_j x_j \right)^2 + \ldots + \left( \sum_j r_{m_n,j} \gamma_j x_j \right)^2. \label{eq_rp12}
\end{eqnarray}
Now $$ E\left( \sum_j r_{1,j} \gamma_j x_j \right)^2 =E\left\{ \sum_j r_{1,j}^2 \gamma_j x_j^2 + \sum_{j\neq j^{\prime}} r_{1,j} r_{1,j^{\prime}} \gamma_j \gamma_{j^{\prime}} x_j x_{j^{\prime}} \right\}= 2 \sum_j \gamma_j x_j^2 =2 \|{\bf x}_{\gamma}\|^2, $$
as $E(r_{i,j}^2)=1$ and $E(r_{i,j} r_{i,j^{\prime}})=0$ as $i=1,2,\ldots, m_n$, $j,j^{\prime}=1,2, \ldots, p_n$, and $j \neq j^{\prime}$. 

\vskip15pt
\noindent {\it Proof of part b.} From (\ref{eq_rp12}) we have
\begin{eqnarray}
 var\left( \|\Rg {\bf x}\|^2|\bgm \right) &=& var \left\{ \sum_i \left( \sum_j r_{i,j} \gamma_j x_j \right)^2   \right\}  = \sum_i var \left( \sum_j r_{i,j} \gamma_j x_j \right)^2  \notag \\
 && \hskip20pt  + \sum_{i\neq i^{\prime}} cov \left\{ \left( \sum_j r_{i,j} \gamma_j x_j \right)^2, \left( \sum_j r_{i^{\prime},j} \gamma_j x_j \right)^2\right\}.  \label{eq_rp13}  
\end{eqnarray}
We will consider each term of (\ref{eq_rp13}) one by one. Consider the first term. Note that 
\begin{eqnarray}
var \left( \sum_j r_{i,j} \gamma_j x_j \right)^2 &=& var \left\{ \sum_j r_{i,j}^2 \gamma_j x_j^2 + \sum_{j\neq k} r_{i,j} r_{i,j^{\prime}} \gamma_j \gamma_k x_j x_{j^{\prime}} \right\}  \notag \\
&=&  var \left\{ \sum_j r_{i,j}^2 \gamma_j x_j^2 \right\} + var \left\{ \sum_{j\neq j^{\prime}} r_{i,j} r_{i,j^{\prime}} \gamma_j \gamma_k x_j x_{j^{\prime}} \right\} \notag \\
&& \hskip25pt + cov \left\{ \sum_j r_{i,j}^2 \gamma_j x_j^2 ,  \sum_{j\neq j^{\prime}} r_{i,j} r_{i,j^{\prime}} \gamma_j \gamma_{j^{\prime}} x_j x_{j^{\prime}} \right\}. \notag
\end{eqnarray}
Consider the first term in (\ref{eq_rp13}).
\begin{eqnarray*}
 var \left\{ \sum_j r_{i,j}^2 \gamma_j x_j^2 \right\} &=& \sum_j var \left( r_{i,j}^2 \gamma_j x_j^2 \right) + \sum_{j \neq j^{\prime} } cov \left( r_{i,j}^2 \gamma_j x_j^2, r_{i,j^{\prime}}^2 \gamma_{j^{\prime}} x_{j^{\prime}}^2 \right)\\
 &=& \sum_j \gamma_j x_j^4 var \left( r_{i,j}^2  \right) + \sum_{j \neq j^{\prime} } \gamma_j x_j^2 \gamma_{j^{\prime}} x_{j^{\prime}}^2 cov \left( r_{i,j}^2 , r_{i,j^{\prime}}^2  \right)\\
   &=& \sum_j \gamma_j x_j^4  \left\{ E\left( r_{i,j}^4 \right) - E^{2}  \left( r_{i,j}^2 \right)   \right\} 
 = 2\left( \frac{1}{\psi}-1  \right) \sum_j \gamma_j x_j^4 .
\end{eqnarray*}

Again,
\begin{eqnarray*}
   var \left\{ \sum_{j\neq j^{\prime}} r_{i,j} r_{i,j^{\prime}} \gamma_j \gamma_k x_j x_{j^{\prime}} \right\}   = E\left(\sum_{j\neq j^{\prime}} r_{i,j} r_{i,j^{\prime}} \gamma_j \gamma_k x_j x_{j^{\prime}} \right)^2 ,\hspace{1.2 in}\\
  = \sum_{j\neq j^{\prime}}  \gamma_j \gamma_k x^2_j x^2_{j^{\prime}} E \left( r_{i,j}^2 r_{i,j^{\prime}}^2 \right)+ \sum_{\substack{\text{$(j,j^{\prime})\neq (k,k^{\prime})$} \\ \text{$j\neq j^{\prime}, k \neq k^{\prime}$}}}  \gamma_j \gamma_k \gamma_{j^{\prime}} \gamma_{k^{\prime}} x^2_j x^2_{j^{\prime}} x^2_k x^2_{k^{\prime}} E \left( r_{i,j} r_{i,j^{\prime}} r_{i,k} r_{i,k^{\prime}} \right)\\
  = 4 \sum_{j\neq j^{\prime}}  \gamma_j \gamma_k x^2_j x^2_{j^{\prime}},\hspace{3.6 in}
\end{eqnarray*}
as the other term will be zero. Next
\begin{eqnarray*}
 &cov \left\{ \sum_j r_{i,j}^2 \gamma_j x_j^2 ,  \sum_{j\neq j^{\prime}} r_{i,j} r_{i,j^{\prime}} \gamma_j \gamma_{j^{\prime}} x_j x_{j^{\prime}} \right\} \hspace{1 in} \\
 & = \sum_j \sum_{k\neq k^{\prime}}  \gamma_j x_j^2 ,    \gamma_k \gamma_{k^{\prime}} x_k x_{k^{\prime}} cov\left(r_{i,j}^2 ,r_{i,k} r_{i,k^{\prime}}\right)=0.
\end{eqnarray*}
Therefore the first term in (\ref{eq_rp13}) is
\begin{eqnarray}
 \sum_i var \left( \sum_j r_{i,j} \gamma_j x_j \right)^2 &=& \frac{1}{m_n} \left[ ( n^{\kappa} -1 ) \sum_j \gamma_j x_j^4 +  \sum_{j\neq j^{\prime}}  \gamma_j \gamma_k x^2_j x^2_{j^{\prime}} \right] \notag \\
 &=& 2 \left[ \left( \frac{1}{\psi} -1 \right) \sum_j \gamma_j x_j^4 + 2 \left(\sum_{j}  \gamma_j  x^2_j \right)^2 \right]. \label{eq_rp15}
\end{eqnarray}
Now consider the last term in (\ref{eq_rp13}).
\begin{eqnarray*}
  && cov \left\{ \left( \sum_j r_{i,j} \gamma_j x_j \right)^2, \left( \sum_j r_{i^{\prime},j} \gamma_j x_j \right)^2\right\}\\
  && =cov \left\{ \sum_j r^2_{i,j} \gamma_j x^2_j + \sum_{j \neq j^{\prime}} r_{i,j} r_{i,j^{\prime}} \gamma_j \gamma_{j^{\prime}} x_j^2 x^2_{j^{\prime}} , \sum_k r^2_{i^{\prime},k} \gamma_k x^2_k + \sum_{k \neq k^{\prime}} r_{i^{\prime},k} r_{i^{\prime},k^{\prime}} \gamma_k \gamma_{k^{\prime}} x_k^2 x^2_{k^{\prime}} \right\} \notag \\
  \end{eqnarray*}
 \begin{eqnarray}
 && =cov \left\{ \sum_j r^2_{i,j} \gamma_j x^2_j, \sum_k r^2_{i^{\prime},k} \gamma_k x^2_k \right\} + cov \left\{ \sum_j r^2_{i,j} \gamma_j x^2_j, \sum_{k \neq k^{\prime}} r_{i^{\prime},k} r_{i^{\prime},k^{\prime}} \gamma_k \gamma_{k^{\prime}} x_k^2 x^2_{k^{\prime}} \right\} \notag \\
 &&+ cov \left\{ \sum_{j \neq j^{\prime}} r_{i,j} r_{i,j^{\prime}} \gamma_j \gamma_{j^{\prime}} x_j^2 x^2_{j^{\prime}} , \sum_k r^2_{i^{\prime},k} \gamma_k x^2_k \right\} \notag \\
 && + cov \left\{ \sum_{j \neq j^{\prime}} r_{i,j} r_{i,j^{\prime}} \gamma_j \gamma_{j^{\prime}} x_j^2 x^2_{j^{\prime}} ,  \sum_{k \neq k^{\prime}} r_{i^{\prime},k} r_{i^{\prime},k^{\prime}} \gamma_k \gamma_{k^{\prime}} x_k^2 x^2_{k^{\prime}}  \right\}. \label{eq_rp14}
\end{eqnarray}
Consider the first term of (\ref{eq_rp14}).
\begin{eqnarray*}
 cov \left\{ \sum_j r^2_{i,j} \gamma_j x^2_j, \sum_k r^2_{i^{\prime},k} \gamma_k x^2_k \right\} = \sum_j \sum_k \gamma_j x^2_j \gamma_k x^2_k cov \left(r^2_{i,j}, r^2_{i^{\prime},k} \right) =0 .
\end{eqnarray*}
Similarly, it can be shown that all other terms in (\ref{eq_rp14}) are zero. Combining the above result and (\ref{eq_rp15}) the proof follows.
\end{proof}

\end{result}
%